\input amstex
\input amsppt.sty
\NoBlackBoxes
\nologo
\parskip 5pt minus 2pt
\hoffset.4truein\voffset-.125truein
\hsize5.75truein\vsize8.5truein

\redefine\i{{\imath}}

\redefine\C{{\Bbb C}}
\redefine\F{{\Bbb F}}
\redefine\N{{\Bbb N}}
\redefine\Q{{\Bbb Q}}
\redefine\R{{\Bbb R}}
\redefine\Z{{\Bbb Z}}

\define\br#1{{\left<{#1}\right>}}
\define\edge#1{{\,\overline{\!#1}}}
\define\endpf{{\hfill $\blacksquare$}}
\redefine\gcd{{\text{\rm gcd}}}
\define\id{{\text{\rm id}}}
\define\isim{{\;@>\,_\sim\,>>\;}}
\define\isom{{\;\simeq\;}}
\redefine\ker{{\text{\rm ker}}}
\define\lcm{{\text{\rm lcm}}}
\redefine\mod{{\text{\rm mod\,}}}
\define\w#1{{\widehat{#1}}}

\define\Pf{{\flushpar{\it Proof.\;\;}}}

\flushpar
{\bf
A family of $p$-adic isometries, fixed points, and the number three}
\newline
Eric S. Brussel
\footnote{
Department of Mathematics and Computer Science,
Emory University, Atlanta, GA, 30322
\newline $\hphantom{\quad}$
The author was partially supported by NSA Grant Number H98230-05-1-0248.
}

\title
{\bf $p$-adic Isometries}
\endtitle

\author
Brussel
\endauthor

\

\flushpar
{\bf Abstract:}
We study the $p$-adic interpolation $\iota_q$ of the arithmetic function
$n\mapsto 1+q+\cdots+q^{n-1}$, where $q\equiv 1(\mod p)$.
We show $\i_q$ has a nontrivial $p$-adic fixed point $z_q$
if and only if $p=3$, $q\not\equiv 1(\mod 9)$, and $q$ is not equal to one of two
$3$-adic integers, $q_0$ and $q_1$.
Setting $\Phi(q)=z_q$, $\Phi(q_0)=0$, and $\Phi(q_1)=1$, we obtain a homeomorphism
$\Phi:U^{(1)}-U^{(2)}\to 3\Z_3\cup(1+3\Z_3)$.
Underlying $\Phi$ are two isometries of the $3$-adic unit disk, which we conjecture are
rigid analytic.

\

\flushpar

\

\flushpar
{\bf Introduction.}

We start with an example from complex analysis.
Let $D$ be the unit disk in the complex plane $\C$.
An {\it isometry} of $D$ is a continuous, 
distance-preserving map from $D$ to $D$.
All analytic isometries of $D$ are rotations, and preserve the complex norm.
They are parameterized in a natural way by $\R/\Z$,
with $t\in\R/\Z$ corresponding to the rotation $\rho_t:z\mapsto ze^{2\pi i t}$.  
The quotient topology on $\R/\Z$ makes the isometries
into a {\it continuous family}, since for all $z\in D$,
$\lim_{t\to t_0}\rho_t(z)=\rho_{t_0}(z)$.
The fixed point set of this family is uninteresting,
since a nontrivial rotation fixes only the origin.
A more interesting set of fixed points is provided by
the larger family of {\it analytic automorphisms} of $D$.
By Schwarz's Lemma, this family is continuously parameterized by $\R/\Z\times D$,
with $(t,z_0)$ corresponding to the mobius transformation
$z\mapsto e^{2\pi it}\frac{z_0-z}{1-\bar z_0 z}$.
A direct computation shows
that an analytic automorphism has either one interior fixed point,
or one boundary fixed point, or two boundary fixed points.

In this paper we study isometries of the {\it $p$-adic} 
unit disk $\Z_p$, and their fixed points.
Let $p$ be a prime, and let $\Z_p$ denote the additive group of $p$-adic integers.
We consider a continuous family of norm-preserving isometries
$$
\i_q:\Z_p\;\longrightarrow\;\Z_p
$$
parameterized by the elements $q$ of 
the topological group $U^{(1)}=1+p\Z_p$ if
$p$ is odd, and $U^{(2)}=1+4\Z_2$ if $p=2$.
Each $\i_q$ is an interpolation of the arithmetic function
on $\N\cup\{0\}$ given by
$$
\i_q(n)=1+q+q^2+\cdots+q^{n-1}.
$$
It is sometimes called the {\it $q$-analog}, or {\it $q$-extension},
of the identity function, and its values are {\it $q$-numbers}.
It is proved in \cite{C} that $\i_q$ is part of a normal basis for the space of
continuous functions from $\Z_p$ to $\Z_p$, along with the other
$q$-binomial coefficients.
In fact, $\i_q$ is also the canonical topological generator of 
the group of continuous 1-cocycles $Z^1(\Z_p,\Z_p)$,
where the $\Z_p$ action on $\Z_p$ is defined by $z*a=q^z a$.

Our goal is to determine $\i_q$'s fixed points.
The reader can immediately verify that $\i_q(0)=0$ and $\i_q(1)=1$; 
we call these fixed points {\it trivial}.
There exist nontrivial fixed points:
if $p=3$ then $-1/2\in\Z_3$,
and it is easily checked that $\i_4(-1/2)=-1/2$.

{\it Results.}
We prove that if $p\neq 3$ or $q\equiv 1(\mod p^2)$ then $\i_q$ has no nontrivial
fixed points.
However, if $p=3$, $q\equiv 1(\mod 3)$, and $q\not\equiv 1(\mod 9)$,
then $\i_q$ has a unique nontrivial fixed point $z_q\in\Z_3$ for all
$q$, with two exceptions.  
The two exceptions, which we call $q_0$ and $q_1$,
canonically determine the ``trivial'' fixed points $z_{q_0}=0$ and $z_{q_1}=1$,
respectively.

The assignment $q\mapsto z_q$, taking an element $q$ of the parameterizing
space to the unique nontrivial fixed point of $\i_q$,
defines a canonical {\it homeomorphism}
$$
\Phi:U^{(1)}-U^{(2)}\;\longrightarrow\;3\Z_3\cup(1+3\Z_3).
$$
Even though we know $\Phi(4)=-1/2$,
we have no closed form expression for $\Phi(q)$.
However, we can show that
underlying $\Phi$ is a pair of {\it isometries}.
That is, first decomposing 
$U^{(1)}-U^{(2)}=(7+9\Z_3)\cup(4+9\Z_3)$, 
we find $\Phi$ takes
$7+9\Z_3$ onto $3\Z_3$, and $4+9\Z_3$ onto $1+3\Z_3$.
Then we prove the compositions 
$$
\align
G&:\Z_3\isim 7+9\Z_3 @>\;\Phi\;>> 3\Z_3\isim\Z_3
\\
F&:\Z_3\isim 4+9\Z_3 @>\;\Phi\;>> 1+3\Z_3\isim\Z_3
\endalign
$$ 
are isometries.  We conjecture, after a suggestion by Tate, that
these functions are {\it rigid analytic}.

In determining the $p$-adic fixed points
we simultaneously determine the {\it modular} fixed points of the induced maps
$$
[\i_q]_{p^n}:\Z_p\;\longrightarrow\;\Z/p^n\Z
$$
for various $n$, defined to be those elements $z\in\Z_p$ such that 
$\i_q(z)\equiv z(\mod p^n)$.  
The cocycle $[\i_q]_{p^n}$ comes up frequently in applications.
For example, if $G$ is a group, $\mu_{p^n}$ is a (multiplicative) $G$-module
of exponent $p^n$, $f$ is a 1-cocycle with values in $\mu_{p^n}$,
and $s\in G$ acts on $\mu_{p^n}$ as exponentiation-by-$q$,
then $f(s^z)=f(s)^{[\i_q(z)]_{p^n}}$ for all $z\in\Z$.
This situation arises over finite fields $\F_q$ that
contain $p$-th roots of unity.

We easily deduce the modular fixed points in all cases except when
$p=3$ and $q\in U^{(1)}-U^{(2)}$, that is, in all cases except when
the isometry $\i_q$ belongs to the family of isometries
that possess nontrivial $p$-adic fixed points. 
These modular fixed points exhibit a regular pattern.
However, when $p=3$ and $q\in U^{(1)}-U^{(2)}$,
the fixed points of $[\i_q]_{3^n}$ exhibit a remarkable, seemingly
erratic pattern that turns out to be
governed completely by the canonical $3$-adic fixed point $z_q=\Phi(q)$.
For example, if $v_0$ is the exponent of the largest power of $3$ to divide $z_q$
or $z_q-1$, and $n$ is sufficiently large,
then the residue of $z_q$ modulo $3^{n-v_0-1}$ is a fixed point
for $[\i_q]_{3^n}$.

Finally, we count the number of modular fixed points of $[\i_q]_{p^n}$
for all primes $p$.  
When $p\neq 3$ or $q\not\in\{q_0,q_1\}$, the number is a certain constant
(which we compute) for all sufficiently large $n$.
For $p=3$ and $q\in\{q_0,q_1\}$, the number of fixed points
for $[\i_q]_{3^n}$ grows without bound as $n$ goes to infinity.

Isometries on $\Z_p$ or on locally compact connected
one-dimensional abelian groups are studied in
\cite{A}, \cite{B}, and \cite{Su}.
More generally, 
investigations into the structure of the space of continuous functions
$C(K,\Q_p)$, where $K$ is a local field, is part of $p$-adic analysis,
and was initiated by Dieudonn\'e in \cite{D}.
Mahler constructed an explicit basis for this space in \cite{M}. 
The concept of $q$-numbers seem to have originated with Jackson, see \cite{J},
and has spawned an industry.
In \cite{F} Fray proved $q$-analogs of theorems of Legendre, Kummer, and
Lucas on $q$-binomial coefficients.  
In \cite{C}
Conrad proved that the set of all $q$-binomial coefficients form a basis
for $C(\Z_p,\Z_p)$.

\

\flushpar
{\bf 1. Background and Notation.}

If $G$ is a group and $g\in G$, we let $o(g)$ denote the order of $g$ in $G$. 
If $M$ is a $G$-module, we write $Z^1(G,M)$ for the group
of 1-cocycles on $G$ with values in $M$, which are functions $f:G\to M$
satisfying the {\it cocycle condition} $f(st)=f(t)^s + f(s)$.
If $G$ is a topological group and $M$ a topological $G$-module,
we assume our 1-cocycles are {\it continuous}.

Let $p$ be a prime, and let $\Z_p$ denote the additive group of $p$-adic
integers, with additive valuation $v_p$.
If $q\in\Z_p^\times$, the group of $p$-adic units,
let $[q]_{p^n}$ denote the image of $q$ in $(\Z/p^n\Z)^\times$.
Then $o([q]_{p^n})$ is the (multiplicative) order of $q$ in $(\Z/p^n\Z)^\times$.
When $q$ is understood we will frequently set 
$o_{p^n}=o([q]_{p^n})$ to save space.
There is a canonical decomposition
$$
\Z_p^\times = \mu'\times U^{(1)}
$$
where $U^{(1)}=\{u\in\Z_p^\times: u\equiv 1(\mod p)\}$ is the group of
{\it principal units}, and
$\mu'$ is the group of prime-to-$p$-order roots of unity, which
is cyclic of order $p-1$.  Let $\mu$ denote the group of all roots of
unity in $\Z_p^\times$; then if $p$ is odd, $\mu=\mu'$, and if $p=2$,
$\mu=\{\pm 1\}$.
Set
$$
U^{(m)}=\{u\in\Z_p:u\equiv 1(\mod p^m)\}.
$$
We summarize some standard facts about these groups (\cite{S}).
$U^{(1)}$ is a (multiplicative) topological group, with the subspace topology. 
The subgroups $U^{(m)}$ of $U^{(1)}$
form a basis of open neighborhoods of the identity.
$U^{(1)}$ has a canonical continuous $\Z_p$-module structure given by $z*q=q^z$, 
where if $z=\lim_{n\to\infty}z_n$ then $q^z:=\lim_{n\to\infty}q^{z_n}$.
If $p$ is odd then $U^{(1)}$ is a free $\Z_p$-module of rank one, i.e.,
a torsion-free procyclic $\Z_p$-module.
If $p=2$ then $U^{(2)}$ is a torsion-free procyclic $\Z_2$-module,
and there is an isomorphism $U^{(1)}\isom\{\pm 1\}\times U^{(2)}$,
given by $q\mapsto (1,q)$ if $q\in U^{(2)}$,
and $q\mapsto(-1,-q)$ if $q\in U^{(1)}-U^{(2)}$.
In particular, $U^{(m)}$ is procyclic if $p$ is odd and $m\geq 1$, 
or if $p=2$ and $m\geq 2$.
If $m<n$, the quotient $U^{(m)}/U^{(n)}$ is represented by the set
$\{1+a_m p^m+\cdots +a_{n-1}p^{n-1}:0\leq a_i<p\}$.
In particular, $|U^{(m)}/U^{(n)}|=p^{n-m}$.
Thus if $p$ is odd, or if $p=2$ and $m\geq 2$,
$U^{(m)}$ is topologically generated
by the elements of $U^{(m)}-U^{(m+1)}$.

We will use the following theorem that goes back to Legendre.
The proof is not hard (\cite{G}).
If $a\in\N$ has $p$-adic expansion $a=a_0+a_1 p+\cdots+a_r p^r$,
let $s_p(a)=a_0+a_1+\cdots+a_r$, the sum of $a$'s digits in base $p$.
Then 
$$
v_p(a!)=\frac{a-s_p(a)}{p-1}.
$$
Using this formula it is not hard to derive the $p$-value of binomial coefficients
(\cite{G}):
if $b\leq a\in\N$, then
$$
v_p(\tbinom ab)=\frac{s_p(b)+s_p(a-b)-s_p(a)}{p-1}.
$$
We will call this expression {\it Kummer's formula}.
It follows immediately by Kummer's formula that $v_p(\binom{p^n}j)=n-v_p(j)$.
We will  also need to use the Binomial Theorem applied to $p$-adic integers.
If $q=1+X\in U^{(1)}$ and $z\in\Z_p$,
then the $p$-adic integer $\binom zi$ is defined as follows.  
If $z=\lim_{n\to\infty}z_n$, where $z_n$ is the residue of $z(\mod p^n)$,
and $i\in\N$, then
$\tbinom zi:=\lim_{n\to\infty}\tbinom{z_n}i$.
The binomial expansion takes the form
$$
(1+X)^z=\sum_{i=0}^\infty\tbinom zi X^i,
$$
where we set $\binom zi=0$ if $z\in\N$ and $z<i$ (see, e.g., \cite{N, Section 5}).

We set up the proper context for our investigation.
Fix $q\in\Z_p^\times$.
The map
$$
\align
\Z\;\times\;\Z_p\;&\longrightarrow\;\Z_p\\
(m,a)\;&\longmapsto\;q^m a
\endalign
$$
defines a nontrivial action of additive groups.
Let $\i_q\in Z^1(\Z,\Z_p)$ be the canonical 1-cocycle,
defined by $\i_q(1)=1$.
It is easy to see that $\i_q$ generates $Z^1(\Z,\Z_p)$, though
we do not need this fact.
The cocycle condition takes the form 
$\i_q(m+n)=q^m \i_q(n)+\i_q(m)$.
In particular, $\i_q(m)=q\i_q(m-1)+1$, and by induction
we have for all $m\in\N$ the formula
$$
\i_q(m)=1+q+\cdots+q^{m-1}.
$$
For any $n\in\N$ we have an induced action 
$\Z\times\Z/p^n\Z\longrightarrow\Z/p^n\Z$, and a canonical cocycle
$[\i_q]_{p^n}$ with image 
$[\i_q(m)]_{p^n}=[1+q+\cdots+q^{m-1}]_{p^n}$
for all $m\in\N$.
It turns out that for a proper analysis,
we must replace $\Z$ with a procyclic group, by
first interpolating the action from $\Z$ to the profinite completion $\w\Z$, 
and then dividing out by the kernel of the extended action.
We call the resulting procyclic group $C_q$,
and view our canonical cocycle $\i_q$ as an element of $Z^1(C_q,\Z_p)$.
In the case of primary interest in this paper, $C_q=\Z_p$.

To derive $C_q$, we start with $\w\Z$.
Every procyclic group is a quotient of $\w\Z$,
hence any procyclic action on $\Z_p$ is the factorization of a continuous
homomorphism $\varphi:\w\Z\to\text{Aut}(\Z_p)\isom\Z_p^\times$.

\proclaim{Lemma 1.1}
Let $\varphi:\w\Z\to\Z_p^\times$ be a continuous homomorphism, and suppose 
$\varphi(1)=q$.  Let $C_q=\w\Z/\ker(\varphi)$ and $o_p=o([q]_p)$.
Then
$$
C_q\isom\cases
\Z/2\Z &\text{if $p=2$ and $q\in\mu$}\\
\Z/o_p\Z &\text{if $p$ is odd and $q\in\mu$}\\
\Z/o_p\Z\times\Z_p &\text{if $q\not\in\mu$}\\
\endcases
$$
\endproclaim

\Pf
For each $n\in\N$ we have an induced map 
$\varphi_n:\w\Z\to\Z_p^\times/(\Z_p^\times)^n$,
whose kernel is an open subgroup $r_n\w\Z$ for some $r_n\in\N$.
Since an element of $\Z_p^\times$ is $1$ if and only if its image in each
$\Z_p^\times/(\Z_p^\times)^n$ is $1$,
$\ker(\varphi)=\bigcap_n\ker(\varphi_n)=\bigcap_n r_n\w\Z$.
Write $q=\omega u$, where $\omega\in\mu$ and $u$ belongs to the torsion-free
$\Z_p$-module $U^{(1)}$ if $p$ is odd, $U^{(2)}$ if $p=2$.
We have $r\in\ker(\varphi_n)$ if and only if $\omega^r=1$ and $u^r=1$ modulo 
$(\Z_p^\times)^n$.
If $n$ divides $n'$ then $\ker(\varphi_n)$ contains $\ker(\varphi_{n'})$,
therefore we may assume $n$ is divisible by $o(\omega)=e_p$,
where $e_2=2$ and $e_p=o_p$ if $p$ is odd.
Since $U^{(1)}$ is prime-to-$p$ divisible, 
$(U^{(1)})^n=(U^{(1)})^{p^{v_p(n)}}$.
Therefore $r_n=\lcm[e_p,o([u]_{p^{v_p(n)}})]$.
If $u=1$, i.e., $q\in\mu$,
this shows $\ker(\varphi)=e_p\w\Z$, hence $C_q\isom\Z/e_p\Z$.
If $q\not\in\mu$, then $o([q]_{p^{v_p(n)}})$ is a power of $p$ that grows without
bound, hence $\ker(\varphi)=r\w\Z$, 
where $r=\lim_{n\to\infty}\lcm[e_p,p^n]$. 
Hence if $p=2$, $C_q\isom\Z_2=\Z/o_2\Z\times\Z_2$; if $p$ is odd,
$C_q\isom\Z/o_p\Z\times\Z_p$.
This completes the proof.

\endpf

Thus $q\in\Z_p^\times$ gives $\Z_p$ a canonical continuous 
$C_q$-module structure,
$$
\align
C_q\times\Z_p\;&\longrightarrow\;\Z_p\\
(z,a)\;&\longmapsto\;q^z a
\endalign
$$
where $q^z a:=\varphi(z)(a)$.  
Note $\Z$ maps into $C_q$ as $m\mapsto[m]$ in $\Z/2\Z$ or $\Z/o_p\Z$ if $q\in\mu$,
and the meaning of $q^z$ is obvious.
If $q\not\in\mu$ then $\Z$ embeds as $m\mapsto([m],m)\in\Z/o_p\Z\times\Z_p$.
To interpret $q^z$ in this case 
we write $q=\omega u$ with $\omega\in\mu'$ and $u\in U^{(1)}$,
and if $z=([a],b)\in \Z/o_p\Z\times\Z_p$, we have $q^z=\omega^a u^b$.

We now replace our original function, which was defined on $\Z$,
with the canonical cocycle $\i_q\in Z^1(C_q,\Z_p)$,
and identify $C_q$ with the groups listed in Lemma 1.1.
Because of the general way in which we constructed $C_q$, 
doing so does not compromise any of the function's properties.
Note $C_q=\Z_p$ if $p$ is odd and $q\in U^{(1)}$, or if $p=2$ and $q\in U^{(2)}$,
and then we have an obvious notion of fixed point.
By replacing $\Z$ by $\Z_p$ in this case,
we have created a more natural setting in which to consider fixed points,
and this proves to be a crucial step for the theory.

We will often cite the following easy observations, 
which we make into a lemma.

\proclaim
{Lemma 1.2}
For all $z\in\N$, we have
$\i_q(0)=0$ and 
$\i_q(-z)=-q^{-z}\i_q(z)=-q^{-1}\i_{q^{-1}}(z)$.
For all $r,z\in C_q$,
$(q-1)\i_q(z)=q^z-1$ and
$\i_q(rz)=\i_{q^z}(r)\i_q(z)$.
\endproclaim

\Pf
The 1-cocycle condition 
yields $\i_q(0)=\i_q(0+0)=\i_q(0)+\i_q(0)$, so $\i_q(0)=0$.
It follows that $0=\i_q(z-z)=q^z \i_q(-z)+\i_q(z)$,
so $\i_q(-z)=-q^{-z}\i_q(z)$.
If $z\in\N$ then
$$
q^{-z}\i_q(z)= q^{-1}+\cdots+q^{-z+1}+q^{-z}=
q^{-1}(1+q^{-1}+\cdots+(q^{-1})^{z-1})=q^{-1}\i_{q^{-1}}(z).
$$
To show $(q-1)\i_q(z)=q^z-1$ for all numbers $z\in\N$ is elementary,
and repeated application of the 1-cocycle condition proves 
$\i_q(rz)=\i_{q^z}(r)\i_q(z)$ for integers $r,z\in\N$.
Then since $\N$ is dense in $C_q$ and the relevant functions
$\i_q$, $q\mapsto q^z$, multiplication, and addition
are continuous, the results extend to $C_q$.

\endpf

\

\flushpar
{\bf 2. Kernel of $\i_q$ and $[\i_q]_{p^n}$.}

We will show that $\i_q$ is injective.
The kernel of $[\i_q]_{p^n}$ is an open (and closed) subset
since $[\i_q]_{p^n}$ is continuous and its image is finite.
It is also a subgroup:  if $z,z'\in\ker([\i_q]_{p^n})$ then
$\i_q(z+z')\equiv q^z\i_q(z')+\i_q(z)\equiv 0(\mod p^n)$, 
and by Lemma 1.2, $\i_q(-z)\equiv-q^{-z}\i_q(z)\equiv 0(\mod p^n)$.  
The formula $\i_q(z)=\i_q(z')+q^z\i_q(z'-z)$ shows that
$\i_q(z)\equiv\i_q(z')(\mod p^n)$ if and only if $z-z'\in\ker([\i_q]_{p^n})$,
so $[\i_q]_{p^n}$ is injective on the quotient $\edge C_q:=C_q/\ker(\i_q)$.
To compute $\edge C_q$ we need a couple of elementary results
on the multiplicative order of $q$.

\proclaim
{Definition 2.1}
{\rm
Let $q\in\Z_p^\times$.
Set $m_0=v_p(q^{o_p}-1)$, and $l_0=v_2(q+1)$, where
$o_p=o([q]_p)$, the multiplicative order of $q$ modulo $p$.
}
\endproclaim
Note $m_0=\infty$ if and only if $p$ is odd and $q\in\mu$,
and $l_0=\infty$, if and only if $p=2$ and $q\in\mu$.
We identify $p^{\infty}$ with $0$.
A quick check shows that if $m_0\neq\infty$ then 
$q^{o_p}\in U^{(m_0)}-U^{(m_0+1)}$.

\proclaim
{Lemma 2.2}
Suppose $q\in\Z_p^\times$ and $n\geq 1$.
Then $m_0\geq n$ if and only if $o([q]_{p^n})=o_p$.
If $m_0<n$ then
$$
o([q]_{p^n})=
\cases
o_p\cdot p^{n-m_0}
&\text{if $p$ is odd, or if $p=2$ and $q\in U^{(2)}$}\\
2^{n-l_0}
&\text{if $p=2$, $q\in U^{(1)}-U^{(2)}$, and $q\not\equiv -1(\mod 2^n)$}\\
2
&\text{if $p=2$ and $q\equiv -1(\mod 2^n)$.}\\
\endcases
$$
\endproclaim

\Pf
We have $o([q]_{p^n})=o_p$ if and only if 
$q^{o_p}\equiv 1(\mod p^n)$, i.e., $m_0\geq n$.
Suppose $m_0<n$, so $o([q]_{p^n})> o_p\geq 1$.
If $p$ is odd, or $p=2$ and $q\in U^{(2)}$,
then $U^{(m_0)}/U^{(n)}$ is cyclic, generated by $q^{o_p}$, 
so $o([q^{o_p}]_{p^n})=|U^{(m_0)}/U^{(n)}|=p^{n-m_0}$,
therefore $o([q]_{p^n})=o_p\cdot p^{n-m_0}$.
If $p=2$ and $q\in U^{(1)}-U^{(2)}$, then $o([q]_{2^n})=2\cdot o([q^2]_{2^n})$,
and since $q^2\in U^{(2)}$, this is $2^{n+1-v_2(q^2-1)}$ by the first case.
Since $v_2(q^2-1)=v_2(q-1)+v_2(q+1)=1+l_0$, $o([q]_{2^n})=2^{n-l_0}$.
If $q\equiv -1(\mod 2^n)$ then obviously $o([q]_{2^n})=2$.

\endpf

\proclaim{Proposition 2.3}
Suppose $q\in\Z_p^\times$.
Fix $n\geq 1$.
Then $\i_q:C_q\longrightarrow\Z_p$ is injective, and
$[\i_q]_{p^n}:C_q\longrightarrow\Z/p^n\Z$ is injective 
if and only if $q\in\mu$.
If $q\not\in\mu$, then $[\i_q]_{p^n}$ is injective on the quotient 
$\edge C_q$, where
$$
\edge C_q\isom\cases
\Z/p^n\Z &\text{if $p$ is odd and $q\in U^{(1)}$, if $p=2$ and $q\in U^{(2)}$,
or if $q\in U^{(n)}$}\\
\Z/2 o([q]_{2^n})\cdot\Z &\text{if $p=2$, $q\in U^{(1)}-U^{(2)}$, 
$q\not\equiv -1(\mod 2^n)$, and $q\not\in U^{(n)}$}\\
\Z/o([q]_{p^n})\cdot\Z &\text{if $q\not\in U^{(1)}$,
or if $p=2$, $q\in U^{(1)}-U^{(2)}$, $q\equiv-1(\mod 2^n)$, and $q\not\in U^{(n)}$.}\\
\endcases
$$
\endproclaim

\Pf
If $q=1$ then $C_q=\{0\}$, so $\i_q$ is injective.
If $q\neq 1$ and $\i_q(z)=0$, then by Lemma 1.2, $q^z=1$,
and since $C_q$ acts faithfully on $\Z_p^\times$, 
$z=0$.  Thus $\i_q$ is injective for all $q\in\Z_p^\times$.

For the modular case,
let $m=v_p(q-1)$.
Note that $\i_q(z)\equiv 0(\mod p^n)$
is equivalent to $q^z-1\equiv 0(\mod (q-1)p^n)$, or, 
since $(q-1)p^n\Z_p=p^{m+n}\Z_p$, to $q^z\equiv 1(\mod p^{m+n})$.
Therefore $\i_q(z)\equiv 0(\mod p^n)$ if and only if $z\in o([q]_{p^{m+n}})\cdot C_q$.
Thus 
$$
\ker([\i_q]_{p^n})=o([q]_{p^{m+n}})\cdot C_q.
$$
If $q\in\mu$ then $|C_q|$ divides $\lcm[2,o_p]$, which divides $o([q]_{p^{m+n}})$,
so $\ker([\i_q]_{p^n})=\{0\}$, and $[\i_q]_{p^n}$ is injective on
$C_q$, as claimed.

If $q\not\in\mu$, then $C_q$ is infinite, so $C_q\neq\edge C_q$.
We have already seen that $[\i_q]_{p^n}$ is injective on $\edge C_q$,
so it remains to compute $o([q]_{p^{m+n}})$ using Lemma 2.2 for the various types
of $q$, and to put this together with the definition of $C_q$.

Suppose $q\not\in\mu$ and $q\not\in U^{(1)}$.
Set $o_{p^a}=o([q]_{p^a})$.
Since $q\not\in U^{(1)}$, $m=0$, so $o_{p^{m+n}}=o_{p^n}$.
By Lemma 1.1, $C_q=\Z/o_p\Z\times\Z_p$.
Since $o_p$ divides $o_{p^n}$, $o_{p^{m+n}}C_q$ equals
$o_{p^n}\Z_p$, and since $\gcd(o_p,p)=1$, $o_{p^n}\Z_p=(o_{p^n}/o_p)\Z_p$.
Therefore $\edge C_q\isom\Z/o_p\Z\times\Z/\frac{o_{p^n}}{o_p}\Z\isom\Z/o_{p^n}\Z$,
as desired.

For the rest of the proof we have $q\not\in\mu$ and $q\in U^{(1)}$.
By Lemma 1.1, $C_q=\Z_p$.
We quickly dispense with the $q\in U^{(n)}$ case:
If $q\in U^{(n)}$ then $\i_q(z)\equiv z(\mod p^n)$ for all $z\in\Z_p$,
hence $\edge C_q\isom\Z/p^n\Z$, as desired.
Assume $q\in U^{(1)}-U^{(n)}$.  We {\it claim}
$$
o_{p^{m+n}}=
\cases 
p^n &\text{if $p$ is odd and $q\in U^{(1)}$, or if $p=2$
and $q\in U^{(2)}$}\\
2\cdot o_{2^n}
&\text{if $p=2$, $q\in U^{(1)}-U^{(2)}$, and $q\not\equiv -1(\mod 2^n)$}\\
o_{2^n} &\text{if $p=2$, $q\in U^{(1)}-U^{(2)}$, and $q\equiv -1(\mod 2^n)$}\\
\endcases
$$
This is immediate from Lemma 2.2; we go through it for the reader's
convenience.
Since $q\not\in U^{(n)}$, $m_0=m<n$, and the second part of Lemma 2.2 applies.  
If $p$ is odd, or if $p=2$ and $q\in U^{(2)}$, then 
by Lemma 2.2, $o_{p^{m+n}}=p^{m+n-m}=p^n$, as desired.
Suppose $p=2$ and $q\in U^{(1)}-U^{(2)}$.  If $q\not\equiv -1(\mod 2^n)$,
then $q\not\equiv -1(\mod 2^{m+n})$, so by Lemma 2.2,
$o_{2^{m+n}}=2^{m+n-l_0}=2^m\cdot o_{2^n}$.
Since $q\in U^{(1)}-U^{(2)}$ we have $m=1$, as desired.
Assume $q\equiv -1(\mod 2^n)$.  Clearly $o_{2^n}=2$ and $m=1$.
If $q\equiv -1(\mod 2^{m+n})$ then 
$o_{2^{m+n}}=2$, so $o_{2^{m+n}}=o_{2^n}$, as desired.
If $q\not\equiv -1(\mod 2^{m+n})$ then by Lemma 2.2, $o_{2^{m+n}}=2^{m+n-l_0}$.
Since $q\equiv -1(\mod 2^n)$ and $q\not\equiv -1(\mod 2^{n+1})$,
we compute $l_0=n$, and we obtain 
$o_{2^{m+n}}=2^m=2=o_{2^n}$, as desired.
This proves the claim.

If $p$ is odd and $q\in U^{(1)}$, or if $p=2$ and $q\in U^{(2)}$,
then by the claim, $\edge C_q=\Z_p/p^n\Z_p=\Z/p^n\Z$.
Similarly,
if $p=2$, $q\in U^{(1)}-U^{(2)}$, and $q\not\equiv  -1(\mod 2^n)$, then 
$\edge C_q=\Z/2 o_{2^n}\Z$, and
if $p=2$, $q\in U^{(1)}-U^{(2)}$, and $q\equiv -1(\mod 2^n)$, then
$\edge C_q=\Z/o_{2^n}\Z$.
This completes the proof.

\endpf

\

\flushpar
{\bf 3. Image of $\i_q$.}

We compute $\i_q(C_q)$ and $\i_q(C_q)(\mod p^n)$,
and determine when $\i_q$ is an isometry.

\proclaim
{Theorem 3.1}
Suppose $q\in\Z_p^\times$.
Then $\i_q$ is surjective if and only if $q\in U^{(1)}$ and $p$ is odd,
or $q\in U^{(2)}$ and $p=2$.  Let $q=\omega u$ be the canonical decomposition
of $q$, where $\omega\in\mu'$ and $u\in U^{(1)}$.
The image of $\i_q$ in $\Z_p$ is
$$
\i_q(C_q)=
\cases
\Z_p &\text{if $p$ is odd and $q\in U^{(1)}$, or if $p=2$ and $q\in U^{(2)}$}\\
2^{l_0}\Z_2\cup(1+2^{l_0}\Z_2) &\text{if $p=2$ and $q\in U^{(1)}-U^{(2)}$}\\
\{\tfrac{\omega^i-1}{q-1}\}_{i=0}^{o_p-1}+p^{m_0}\Z_p 
&\text{if $p$ is odd and $q\not\in U^{(1)}$.}\\
\endcases
$$
\endproclaim

\Pf
Suppose $q\in U^{(1)}$ and $p$ is odd, or $q\in U^{(2)}$ and $p=2$.
Then $C_q=\Z_p$ by Lemma 1.1, and $\Z_p*q=U^{(m_0)}$.
Therefore the image of the function $q^x-1:z\mapsto q^z-1$ 
is $p^{m_0}\Z_p=(q-1)\Z_p$.
Since $\i_q(z)=(q^z-1)/(q-1)$, $\i_q(\Z_p)=\Z_p$.

If $p=2$ and $q\in U^{(1)}-U^{(2)}$, then $C_q=\Z_2$ and
$q^2\in U^{(2)}$, and since $v_2(q^2-1)=l_0+1$, $2\Z_2*q=U^{(l_0+1)}$.
Therefore $\Z_2*q=2\Z_2*q\cup(1+2\Z_2)*q=U^{(l_0+1)}\cup qU^{(l_0+1)}$.
Since $q\in\Z_2^\times$, $q U^{(l_0+1)}=q+2^{l_0+1}\Z_2$.
Therefore the image of $q^x-1$ is $2^{l_0+1}\Z_2\cup(q-1+2^{l_0+1}\Z_2)$.
Since $q-1$ is $2$ times a unit in $\Z_2$, we conclude
the image of $\i_q$ is $2^{l_0}\Z_2\cup(1+2^{l_0}\Z_2)$.

If $q\not\in U^{(1)}$, then $C_q*q=\br{\omega}\times U^{(m_0)}
=\br{\omega}+p^{m_0}\Z_p$,
so the image of $q^x-1$ is $(\br{\omega}-1)+p^{m_0}\Z_p$.
Since $q-1$ is a unit, the image of
$\i_q$ is the set $(\br{\omega}-1)/(q-1)+ p^{m_0}\Z_p$,
a finite union of additive cosets.  If $i\not\equiv  j(\mod o_p)$ then 
$\omega^i\not\equiv \omega^j(\mod p)$, so these cosets are all distinct,
and there are exactly $o_p$ of them.
Thus $\i_q$ is surjective if and only if $o_p=[\Z_p:p^{m_0}\Z_p]=p^{m_0}$,
and since $o_p$ is prime to $p$, this proves $\i_q$ is not surjective.

\endpf

We need a technical lemma computing the $p$-value of
the numbers $\i_q(z)$.  The result follows from the $q$-Kummer
theorem proved by Fray in \cite{F}, though this
is not immediately apparent due to the much greater level of generality 
in \cite{F}.

\proclaim
{Lemma 3.2}
Suppose $q\in \Z_p^\times-\{1\}$ and $z\in C_q$.
Then $v_p(\i_q(z))=0$ if and only if $q\in U^{(1)}$ and $z\in\Z_p^\times$,
or $q\not\in U^{(1)}$ and $z\not\in o_p C_q$. 
If $v_p(\i_q(z))\neq 0$, then
$$
v_p(\i_q(z))=\cases
v_p(z')+m_0 &\text{if $p$ is odd and $q\not\in U^{(1)}$}\\
v_p(z) &\text{if $p$ is odd and $q\in U^{(1)}$, or if $p=2$ and $q\in U^{(2)}$}\\
v_2(z)+l_0-1 &\text{if $p=2$ and $q\in U^{(1)}-U^{(2)}$}
\endcases
$$
where $z'\in\Z_p$ is given by $z\mapsto o_p\cdot z'$,
under the canonical isomorphism $o_p C_q\isim\Z_p$.
\endproclaim

\Pf
If $q\in\Z_p^\times-U^{(1)}$, then
$v_p(q^z-1)=v_p(\i_q(z))+v_p(q-1)=v_p(\i_q(z))$, so
$v_p(\i_q(z))=0$ if and only if $q^z\not\in U^{(1)}$, i.e., 
$z\not\in o_p C_q$.
Suppose $q\in U^{(1)}-\{1\}$, so $1\leq m_0=v_p(q-1)<\infty$.
Then $v_p(\i_q(z))=0$ if and only if $v_p(q^z-1)=v_p(q-1)$.
Since $U^{(m)}/U^{(m+1)}$ has order $p$ for all $m\geq 1$,
it is clear that $v_p(q^z-1)=v_p(q-1)$ if and only if $v_p(z)=0$,
i.e., $z\in\Z_p^\times$.

Assume for the rest of the proof that $v_p(\i_q(z))\neq 0$.
Suppose $p$ is odd, then $z\in o_p\cdot C_q$.
Write $z=o_p\cdot z'$, with $z'\in\Z_p$.
Since $q\neq 1$, $q^{o_p}$ generates $U^{(m_0)}$ topologically,
so $v_p(q^z-1)=v_p((q^{o_p})^{z'}-1)=v_p(z')+m_0$,
and $v_p(\i_q(z))=v_p(z')+m_0-v_p(q-1)$.
If $q\in U^{(1)}$ then $m_0-v_p(q-1)=0$ and $v_p(z')=v_p(z)$, 
so $v_p(\i_q(z))=v_p(z)$, as desired.
If $q\in\Z_p^\times-U^{(1)}$, then $v_p(q-1)=0$, so 
$v_p(\i_q(z))=v_p(z')+m_0$, as desired.

Suppose $p=2$, then $z\in 2\Z_2$, and we can write $z=2z'$, with $z'\in\Z_2$.
Since $q\in U^{(1)}$ we have $q^2\in U^{(2)}$, and $q^2$ generates
$U^{(m_0+l_0)}$ topologically.
Thus $v_2(q^z-1)=v_2(z')+l_0+m_0$, so $v_2(\i_q(z))=v_2(z')+l_0=v_2(z)+l_0-1$.
This completes the proof.

\endpf

\proclaim
{Corollary 3.3}
Suppose $q\in\Z_p^\times$.
Then $\i_q$ is an isometry if and only if 
$p$ is odd and $q\in U^{(1)}$, or $p=2$ and $q\in U^{(2)}$.
Every isometry $\i_q$ preserves the norm.
\endproclaim

\Pf
If $\i_q$ is an isometry then it is surjective by definition,
and by Theorem 3.1, 
either $p$ is odd and $q\in U^{(1)}$, or $p=2$ and $q\in U^{(2)}$.
Conversely, suppose $p$ is odd and $q\in U^{(1)}$, or $p=2$ and $q\in U^{(2)}$.
By Theorem 3.1, $\i_q$ is surjective.
By Proposition 2.3, $\i_q$ is injective,
hence it is bijective.
Since $C_q$ is compact, a continuous bijection on $C_q$ is a homeomorphism,
therefore $\i_q$ is a homeomorphism.
By Lemma 3.2, $v_p(\i_q(z))=v_p(z)$ for all $z\in\Z_p$, i.e.,
$\i_q$ preserves the $p$-adic norm.
This completes the proof.

\endpf

\proclaim{Theorem 3.4}
Suppose $q\in\Z_p^\times$.
Then $[\i_q]_{p^n}$ is surjective if and only if $p$ is odd and $q\in U^{(1)}$,
or if $p=2$ and either $n=1$ or $q\in U^{(2)}$.
The image of $[\i_q]_{p^n}$ in $\Z/p^n\Z$ is
$$
\cases
\Z/p^n\Z
&\text{if $p$ is odd and $q\in U^{(1)}$, or if $p=2$ and $q\in U^{(2)}$}\\
\big(2^{l_0}\Z/2^n\Z\big)\cup\big(1+2^{l_0}\Z/2^n\Z\big)
&\text{if $p=2$ and $q\in U^{(1)}-U^{(2)}$}\\
\bigcup_{z_0=0}^{o_p-1}\big(\i_q(z_0)+p^{m_0}\Z/p^n\Z\big)
&\text{if $p$ is odd and $q\not\in U^{(1)}$.}\\
\endcases
$$
We set $p^m\Z/p^n\Z=0$ if $m>n$.
In the last case, the $\i_q(z_0)$ are
all incongruent modulo $p$.
\endproclaim

\Pf
Except for the last case,
the computation of $\i_q(C_q)(\mod p^n)$ is immediate by Theorem 3.1.
For the last case, suppose $p$ is odd and $q\in\Z_p^\times-U^{(1)}$.
Then $q=\omega u$ where $u\in U^{(m_0)}-U^{(m_0+1)}$,
and it follows immediately that for $z_0=0,1,\dots,o_p-1$,
$q^{z_0}\equiv \omega^{z_0}(\mod p^{m_0})$,
hence $\omega^{z_0}-1\equiv q^{z_0}-1(\mod p^{m_0})$.
Since $q-1$ is a unit, $\frac{\omega^{z_0}-1}{q-1}\equiv \i_q(z_0)(\mod p^{m_0})$.
The $\i_q(z_0)$ are incongruent by the proof of Theorem 3.1, so this proves all
but the first statement.
If $p$ is odd and $q\in U^{(1)}$, or if $p=2$ and $q\in U^{(2)}$, then 
$[\i_q]_{p^n}$ is surjective by Theorem 3.1.
If $p=2$ and $n=1$, then an easy computation shows $[\i_q]_{2^n}$ is surjective.
Conversely,
if $p=2$, $n\geq 2$, and $q\in U^{(1)}-U^{(2)}$, then $l_0\geq 2$, and 
$[\i_q]_{2^n}$
is not surjective, and if $p$ is odd and $q\in\Z_p^\times-U^{(1)}$, then 
$[\i_q]_{p^n}$ is not surjective since $o_p<p$.
This completes the proof.

\endpf

\proclaim{Corollary 3.5}
Suppose $p$ is prime, $n\in\N$, and $q\in U^{(1)}$.
Then
$$
\sum_{z=0}^{p^n-1}\i_q(z)\equiv 
\cases
2^{n-1}(\mod 2^n) &\text{if $p=2$}\\
0(\mod p^n) &\text{if $p$ is odd.}
\endcases
$$
\endproclaim

\Pf
If $p$ is odd and $q\in U^{(1)}$,
or if $p=2$ and $q\in U^{(2)}$, 
then by Theorem 3.4 $[\i_q]_{p^n}$ 
is surjective, and $\sum_{z=0}^{p^n-1}\i_q(z)=\sum_{z=0}^{p^n-1}z
=p^n(p^n-1)/2$.  If $p$ is odd then this expression is congruent to $0$
modulo $p^n$, and if $p=2$ it is congruent to $2^{n-1}$.

Suppose $p=2$ and $q\in U^{(1)}-U^{(2)}$.
By Theorem 3.4, $|[\i_q(C_q)]_{2^n}|=2^{n-l_0+1}$,
and the image of $[\i_q]_{2^n}$ is $2^{l_0}\Z/2^n\Z\cup(1+2^{l_0}\Z/2^n\Z)$.
We need to sum the elements of this set, and then
count each element $2^n/|[\i_q(C_q)]_{2^n}|=2^{l_0-1}$ times.
We represent $2^{l_0}\Z/2^n\Z$ by the numbers
$$
S=\left\{2^{l_0}(a_0+a_1\cdot 2+\cdots+a_{n-l_0-1}2^{n-l_0-1}):a_i\in\{0,1\}
\right\}
$$
Then $\i_q(C_q)\equiv S\cup(1+S)(\mod 2^n)$.  
We divide $S$ into unordered pairs $\{s,t\}$, where
if $s=2^{l_0}(a_0+\cdots+a_{n-l_0-1}2^{n-l_0-1})$ then
$t=2^{l_0}(b_0+\cdots+b_{n-l_0-1}2^{n-l_0-1})$, with $b_i=1-a_i$.
Note $s\neq t$, so $S$ is the disjoint union of the pairs $\{s,t\}$,
and
$s+t=2^{l_0}(1+2+\cdots+2^{n-l_0-1})=2^{l_0}(2^{n-l_0-1}-1)$.
Thus the sum of the elements of $S$ is 
$2^{n-l_0}2^{l_0}(2^{n-l_0-1}-1)=2^n(2^{n-l_0-1}-1)$.
Similarly the sum of the elements of $1+S$ is $2^{n-l_0}+2^n(2^{n-l_0-1}-1)$.
The total sum, multiplied by $2^{l_0-1}$, is
$$
2^{l_0-1}(2^{n-l_0}+2^{n+1}(2^{n-l_0-1}-1))
\equiv 2^{n-1}(\mod 2^n).
$$
This completes the proof.

\endpf

\

\flushpar
{\bf 4. Fixed Points.}

If $q\in U^{(1)}$ then $C_q=\Z_p$.
In this case $\i_q$ has an obvious notion of fixed point.

\proclaim{Definition 4.1}
{\rm
Suppose $q=U^{(1)}$.
We say $z\in\Z_p$ is a {\it $p$-adic fixed point} of 
$\i_q$ if $\i_q(z)=z$, 
and a {\it modular fixed point}
of $[\i_q]_{p^n}$ if $\i_q(z)\equiv z(\mod p^n)$.
}
\endproclaim

For example, $\i_q$ fixes $0$ and $1$.
We call these fixed points {\it trivial}.
It is clear that
$z$ is a $p$-adic fixed point for $\i_q$ if and only if $z$ is a modular fixed point
for $[\i_q]_{p^n}$ for all $n$.
The next result shows there are certain modular fixed points that always appear,
and that in most cases these are the only ones.

\proclaim{Theorem 4.2 ``Modular Fixed Point Pairs''}
Suppose $n\in\N$, $q\in U^{(1)}$, and $z\in\Z_p$.
Then $z$ is a modular fixed point of $[\i_q]_{p^n}$ 
if 
$$
v_p(z(z-1))\geq n-v_p(q-1)+v_p(2).
\tag $*$
$$
If either $p\neq 3$, $q\in U^{(2)}$, $n\leq 2$, or $z\equiv 2(\mod 3)$,
then $z$ is a fixed point if and only if $(*)$ holds.
If $p\neq 3$, $q\in U^{(2)}$, or $n\leq 2$,
the complete set of modular fixed points is $a_0\Z_p\cup(1+a_0\Z_p)$, where
$$
a_0=\cases
o([q]_{p^n}) &\text{if $p$ is odd}\\
2\cdot o([q]_{2^n}) &\text{if $p=2$ and $q\in U^{(2)}$}\\
2^n &\text{if $p=2$ and $q\in U^{(1)}-U^{(2)}$.}
\endcases
$$
\endproclaim

\Pf
We set $m_0=v_p(q-1)$, as in Definition 2.1.
It is easily seen that
if $z=0$ or $1$ then $z$ is a fixed point of $[\i_q]_{p^n}$,
and if $z=2$ then $z$ is a fixed point of $[\i_q]_{p^n}$ 
if and only if $q\equiv  1(\mod p^n)$.  
On the other hand if $z=0$ or $1$ then it is immediate that $(*)$ holds,
and if $z=2$ then $(*)$ holds if and only if $n-m_0\leq 0$,
i.e., $q\equiv 1(\mod p^n)$.
Therefore all of the statements hold in these
cases, and we will assume $z\neq 0,1,2$ in the following.

Let $X=q-1$, then $v_p(X)=m_0\geq 1$.
By Definition 4.1, $z$ is a fixed point of $[\i_q]_{p^n}$
if and only if $\i_q(z)\equiv z(\mod p^n)$, and multiplying both
sides by $X$, we see this is equivalent to $(1+X)^z\equiv 1+zX(\mod p^{m_0+n})$.
By the Binomial Theorem we have
$(1+X)^z=\sum_{i=0}^\infty\binom zi X^i$,
so $z$ is a fixed point for $[\i_q]_{p^n}$
if and only if
$\sum_{i=2}^\infty\binom zi X^i\equiv 0(\mod p^{m_0+n})$.
We factor out $X$ and reduce the modulus to $p^n$, 
and replace $z$ by its residue $(\mod p^n)$.  Then we have that
$z$ is a fixed point for $[\i_q]_{p^n}$ if and only if
$$
\sum_{i=2}^z\tbinom zi X^{i-1}\equiv 0(\mod p^n).
\tag $**$
$$
We immediately dispense with the $n\leq 2$ case.  For 
if $n\leq 2$, then $X^2\equiv 0(\mod p^n)$, so $(**)$ becomes 
$\frac{z(z-1)}2 X\equiv 0(\mod p^n)$, hence $z$ is fixed by $[\i_q]_{p^n}$
if and only if $(*)$ holds.

To analyze the sum in $(**)$ in the remaining cases we will prove a
{\it claim}, that the $i=2$ term of $(**)$
has the smallest $p$-value,
and if $p\neq 3$, $m_0\geq 2$, or $z\equiv 2(\mod 3)$,
then the $i=2$ term has the {\it unique} smallest value.
This claim proves the first two statements.
For it implies
the value of the sum in $(**)$ is at least $v_p(z(z-1)X/2)$,
hence that if $(*)$ holds then $[\i_q]_{p^n}$ fixes $z$,
as desired.  Conversely the claim implies that
if $p\neq 3$, $m_0\geq 2$, or $z\equiv 2(\mod 3)$, then the value of the
sum in $(**)$ is exactly $v_0(z(z-1)X/2)$, so 
if $z$ is a fixed point for $[\i_q]_{p^n}$, then $(*)$ holds.

To prove the claim we compute the difference 
in value between the $i=2$ term and the $i=j\geq 3$ term,
$$
v_p\big(\tfrac{z(z-1)\cdots(z-j+1)\cdot 2\cdot p^{(j-1)m_0}}{z(z-1)j!\cdot p^{m_0}}\big)
=v_p\big(\tfrac{(z-2)\cdots(z-j+1)\cdot 2\cdot p^{(j-2)m_0}}{j!}\big).
$$
The claim holds that this number
is always nonnegative, and is positive when $p\neq 3$, $m_0\geq 2$, or $z\equiv 2(\mod 3)$.
Thus we must show
$$
v_p(j!)\leq v_p((z-2)\cdots(z-j+1)\cdot 2\cdot p^{(j-2)m_0})
$$
for $j\geq 3$,
with strict inequality when $p\neq 3$, $m_0\geq 2$, or $z\equiv 2(\mod 3)$.
We resort to a brute force case analysis.
To weed out most of the cases we will use the following bounds.
By Legendre's Theorem, for the denominator we have
$v_p(j!)\leq\lfloor\tfrac{j-1}{p-1}\rfloor$,
since always $s_p(j)\geq 1$ for $j\neq 0$.
On the other hand, for the numerator we have
$$
v_p((z-2)\cdots(z-j+1)\cdot 2\cdot p^{(j-2)m_0})\geq v_p(2\cdot p^{(j-2)m_0})
=(j-2)m_0+v_p(2).
$$
Therefore to prove the claim it is sufficient, but not necessary, to show
$$
\lfloor\tfrac{j-1}{p-1}\rfloor\leq (j-2)m_0+v_p(2)
$$
for $j\geq 3$, with strict inequality when $p\neq 3$, $m_0\geq 2$, or $z\equiv 2(\mod 3)$.

Suppose $p=3$ and $m_0=1$.
If $j\geq 4$ we have $\lfloor\frac{j-1}{p-1}\rfloor=\lfloor\frac{j-1}2\rfloor<j-2$,
so we have strict inequality for $j\geq 4$.
When $j=3$ we have $v_3(j!)=1$ and 
$$
v_p((z-2)\cdots(z-j+1)\cdot 2\cdot p^{(j-2)m_0})=v_3((z-2)\cdot 2\cdot 3)
=v_3(z-2)+1
$$
and we see the former is always less than or equal to the latter, with strict
inequality if and only if $v_3(z-2)\neq 0$, i.e., $z\equiv 2(\mod 3)$.
This proves the first part of the claim,
and the $p=3$, $m_0=1$, $z\equiv 2(\mod 3)$ case of the second part.
It remains to prove strict inequality when $p\neq 3$ or $m_0\geq 2$.

Assume $m_0\geq 2$.
If $p=2$ then $\lfloor\frac{j-1}{p-1}\rfloor=j-1=(j-2)+1<(j-2)m_0+1$,
so we have strict inequality in this case.
If $p$ is odd then $\lfloor\frac{j-1}{p-1}\rfloor<j-1\leq (j-2)2\leq(j-2)m_0$,
since $j\geq 3$, so the claim is true in this case.
Thus we have proved the claim if $m_0\geq 2$.

Suppose $p\neq 3$ and $m_0=1$.
If $p=2$ and $j=3$, then (for all $m_0$) we have
$v_p(j!)=1<j-1=(j-2)+1\leq(j-2)m_0+v_p(2)$, so the claim is true in this case.
If $p=2$ and $j\geq 4$ then $z\geq 4$, hence $v_2(z-2)=1$ or $v_2(z-3)=1$,
and we have 
$$
v_p(j!)\leq\lfloor\tfrac{j-1}{p-1}\rfloor=j-1<j=(j-2)+2
\leq v_p((z-2)(z-3)\cdots(z-j+1)\cdot 2\cdot p^{(j-2)m_0})
$$
and the claim is proved in this case.
If $p\neq 2,3$ then 
$\lfloor\frac{j-1}{p-1}\rfloor\leq\lfloor\frac{j-1}4\rfloor\leq\frac{j-1}4<j-2$,
proving the claim in these cases.
This completes the proof of the claim, and hence of the first two statements.

We now compute the complete set of modular fixed points for $[\i_q]_{p^n}$
when $p\neq 3$, $q\in U^{(2)}$, or $n\leq 2$.
We will determine $a_0$ such that the set of elements that satisfy $(*)$
has the form $a_0\Z_p\cup(1+a_0\Z_p)$.
Set $v_0=v_p(z(z-1))$ and $o_{p^n}=o([q]_{p^n})$.
If $q\equiv 1(\mod p^n)$ then $(*)$ is satisfied for all $z\in\Z_p$,
and since $o_{p^n}=1$, we have $\Z_p=a_0\Z_p\cup(1+a_0\Z_p)$
if $a_0=o_{p^n}$ if $p$ is odd, or $a_0=2\cdot o_{2^n}$ if $p=2$ and $q\in U^{(2)}$,
or $a_0=2^n=2$ if $q\in U^{(1)}-U^{(2)}$, as claimed.
Therefore we will now assume $q\not\equiv  1(\mod p^n)$.
First we consider the case $p\neq 3$, then $p=3$ and $q\in U^{(2)}$,
and finally $p\neq 3$, $q\in U^{(1)}-U^{(2)}$, and $n\leq 2$.

Suppose $p\neq 3$ and $p$ is odd.  Then $v_p(2)=0$, and $(*)$ becomes
$v_0\geq n-m_0$.  By Lemma 2.2, $m_0=n-v_p(o_{p^n})$,
so this reduces to $v_0\geq v_p(o_{p^n})$, i.e., 
$z\in o_{p^n}\Z_p\cup(1+o_{p^n}\Z_p)$.  Thus we set $a_0=o_{p^n}$ in this case.
Suppose $p=2$ and $q\in U^{(1)}-U^{(2)}$.  Then $m_0=1$ and $v_2(2)=1$,
so $(*)$ becomes $v_0\geq n$.  Therefore $z$ satisfies $(*)$ if and 
only if $z\in 2^n\Z_2$ or $z\in 1+2^n\Z_2$, and we set $a_0=2^n$.
Suppose $p=2$ and $q\in U^{(2)}$.  
Then by Lemma 2.2, $m_0=n-v_2(o_{2^n})$, and $(*)$ becomes
$v_0\geq v_2(o_{2^n})+1=v_2(2\cdot o_{2^n})$.
This holds for $z$ if and only if $z\in 2\cdot o_{2^n}\Z_2$ or $z\in 1+2\cdot o_{2^n}\Z_2$,
so we take $a_0=2\cdot o_{2^n}$.
Suppose $p=3$ and $q\in U^{(2)}$.  Then by Lemma 2.2 we have
$m_0=n-o_{3^n}$, and $(*)$ becomes $v_0\geq v_3(o_{3^n})$,
which holds if and only if $z\in o_{3^n}\Z_3$ or $z\in 1+o_{3^n}\Z_3$.
Thus we take $a_0=o_{3^n}$.
Finally, suppose $p\neq 3$, $q\in U^{(1)}-U^{(2)}$, and $n\leq 2$.  
Since we assume $q\not\equiv  1(\mod p^n)$, we have $n=2$.
Then $(*)$ becomes $v_0\geq n-1=1$, which holds
if and only if $z\in 3\Z_3$ or $z\in 1+3\Z_3$.  Since $3=o_{p^n}$,
we may take $a_0=o_{p^n}$.
This completes the proof.

\endpf

\flushpar
{\bf Remark 4.3.}
It is quickly seen that
the criterion $(*)$ fails to give all fixed points when 
$q\in U^{(1)}-U^{(2)}$ and $p=3$.
For example, if $p=3$, $q=4$, and $n=3$, then $z=4$ is a fixed point of $[\i_q]_{3^n}$
not indicated by the criterion.  For we compute $v_3(4\cdot 3)=1<3-v_3(3)=2$,
so $z=4$ does not satisfy $(*)$, yet
$\i_4(4)\equiv (4^4-1)/(4-1)\equiv 255/3\equiv 85\equiv  4(\mod 3^3)$.

\proclaim{Corollary 4.4}
Suppose $q\in U^{(1)}-\{1\}$.
If $p\neq 3$ or $q\in U^{(2)}$ 
then $\i_q$ has only the trivial $p$-adic 
fixed points $0$ and $1$.
\endproclaim

\Pf
By Theorem 4.2 and the hypotheses,
any $p$-adic fixed point $z$ of $\i_q$ is a fixed point of $[\i_q]_{p^n}$
for all $n$, and so belongs to the set
$a_0\Z_p\cup(1+a_0\Z_p)$ for all $n$, where $a_0$ is as in Theorem 4.2.
If $p$ is odd we have $a_0=o([q]_{p^n})$.
Since $q\in U^{(1)}-\{1\}$, $q$ is not a root of unity,
and in $\Z_p$ we have
$\lim_{n\to\infty}o([q]_{p^n})=0$, and we conclude $z=0$ or $z=1$.
If $p=2$ and $q\in U^{(2)}$ we have $a_0=2o([q]_{2^n})$,
and since again $q$ is not a root of unity, a similar argument holds.
If $p=2$ and $q\in U^{(1)}-U^{(2)}$ we have $a_0=2^n$, and
since $\lim_{n\to\infty}2^n=0$, we again conclude $z=0$ or $z=1$.

\endpf

\

\flushpar
{\bf 5. The Number Three.}

By Corollary 4.4, the only $q\in U^{(1)}-\{1\}$ for
which $\i_q$ could conceivably have a nontrivial $p$-adic fixed point
are those that generate $U^{(1)}$ when $p=3$.
Moreover, it follows by Theorem 4.2 that any $3$-adic fixed point $z$ for
$\i_q$ satisfies $v_3(z(z-1))\geq 1$.
Therefore we must have $q\in U^{(1)}-U^{(2)}$ and $z\in 3\Z_3\cup(1+3\Z_3)$.
We will prove the following result, which combines Lemma 5.15 and Theorem 6.1 below.

\proclaim
{Theorem}
There exist unique elements $q_0\equiv 7(\mod 9)$ and
$q_1\equiv 4(\mod 9)$ in $U^{(1)}-U^{(2)}$ such that
$\i_{q_0}$ and $\i_{q_1}$ have no $3$-adic fixed points.
If $q\not\in\{q_0,q_1\}$, then $\i_q$ has a unique nontrivial
$3$-adic fixed point $z_q$.
If $q\in 4+9\Z_3$ then $z_q\in 1+3\Z_3$; if $q\in 7+9\Z_3$ then $z_q\in 3\Z_3$.
\endproclaim

Instead of trying to construct the fixed points for each 
$[\i_q]_{3^n}$ directly,
and then piecing them together to find $3$-adic fixed points,
our strategy is to first construct a $q$ for each 
$z\in 3\Z_3\cup(1+3\Z_3)-\{0,1\}$, such that $\i_q(z)\equiv z(\mod 3^n)$.  
This is Proposition 5.2.
Passing to the limit, we obtain a unique $q$ for each of these $z$,
such that $\i_q(z)=z$.  This is Corollary 5.4.
In Theorem 5.8 we use these results to establish
the two possible structures of the set of all fixed points for $[\i_q]_{3^n}$,
for any given $q\in U^{(1)}-U^{(2)}$.
Which of these two structures applies depends on whether $\i_q$ has
a ``rooted fixed point'' modulo $3^n$.
By piecing these sets together and applying
some subtle counting arguments, we establish in Theorem 5.13 the existence of a
uniquely determined $3$-adic fixed point for those $\i_q$ that exhibit an
rooted fixed point modulo $3^n$ for some $n$.
This accounts for those $q$ constructed in Corollary 5.4.
To show that these are all of the $q$ in $U^{(1)}-U^{(2)}$, 
with two exceptions in Lemma 5.15,
we resort to topological considerations.
This is Theorem 6.1, the main theorem of the paper.

Set $v=v_3$, the additive valuation on $\Z_3$.
The following lemma provides the technical explanation for the exceptional role of
the number three in our context.
Let $\Z[x_1,x_2,\dots]$ denote the polynomial ring in indeterminates $\{x_i\}_{i\in\N}$.
Extend $v$ to this ring by setting $v(x_i)=0$ for all $i$,
and let $\Z[x_1,x_2,\dots]_3$ denote the resulting completion with respect to $3$.
Let $X=x_1 3+x_2 3^2+\cdots\in\Z[x_1,x_2,\dots]_3$, 
and suppose $z\in\Z_3$ satisfies $1\leq v(z(z-1))<\infty$.
Set
$$
S=S(z)=\sum_{j=2}^\infty\tbinom zj X^j,
$$
where $\binom mj=0$ if $m\in\N$ and $m<j$.
For each $j$ we have a $3$-adic expansion $X^j=\sum_{i=j}^\infty a_i 3^i$,
where $a_i\in\Z[x_1,\dots,x_{i-j+1}]$.
This polynomial ring is free on the monomials in the $x_i$, so $a_i$ is uniquely
expressible as a sum of monomials in $x_1$ through $x_{i-j+1}$, with coefficients
in $\{0,1,2\}$.
Fix $k$. 
For each $j$,
let $l_j$ be the smallest number such that $x_k$ appears in $a_{l_j}$,
and let $b_{j,k}$ be the sum of those monomials in $a_{l_j}$ in which $x_k$ appears.
We call $\binom zj b_{j,k}3^{l_j}$ the {\it minimal $x_k$-term} of $\binom zj X^j$.
It represents the additive factor of $\binom zj X^j$ of smallest homogeneous 
$3$-value that is divisible by $x_k$.
Suppose 
$w_0=\inf_j v(\binom zj b_{j,k}3^{l_j})=\inf_j(v(\binom zj)+l_j)$.
The {\it minimal $x_k$-term in $S$} is the sum of those minimal $x_k$-terms
$\binom zj b_{j,k}3^{l_j}$ such that $v(\binom zj)+l_j=w_0$.
It represents the additive factor of $S$ of smallest homogeneous $3$-value that is
divisible by $x_k$.
Note that since the monomials form a basis, adding minimal $x_k$-terms
of given (minimal) value does not change that value, so the minimal $x_k$-term
of $S$ has value $w_0$.

\proclaim
{Lemma 5.1}
Suppose $z\in\Z_3$ satisfies $1\leq v(z(z-1))<\infty$.
Then the minimal $x_k$-term of $S$ has value $v(z(z-1))+k+1$.
The minimal $x_1$-term is $\binom z2 x_1^2 3^2+\binom z3 x_1^3 3^3$,
and for $k\geq 2$ the minimal $x_k$-term is $\binom z2 2x_1 x_k 3^{k+1}$.
\endproclaim

\Pf
Set $v_0=v(z(z-1))$.
We treat the $k=1$ case first.
Let $Y=X-x_1 3\in \Z[x_1,x_2,\dots]_3$, then $v(Y)=2$.
By inspection, the term of smallest $3$-value in 
$X^2=(Y+x_1 3)^2=Y^2+2x_1 3Y+x_1^2 3^2$ is $x_1^2 3^2$.
Therefore the minimal $x_1$-term in 
$\binom z2 X^2$ is $\binom z2 x_1^2 3^2$.
If $j=3$ then similarly the minimal $x_1$-term
in $\binom z3 X^3=\binom z3(Y+x_1 3)^3$ is $\binom z3 x_1^3 3^3$.
Since $v(\binom z2)=v_0$ and $v(\binom z3)=v_0-1$,
$v(\binom z2 x_1^2 3^2)=v(\binom z3 x_1^3 3^3)=v_0+2$.

To prove the first statement we must show
the minimal $x_1$-terms of $\binom zj X^j$ for $j\geq 4$ have higher $3$-value.
It is easy to see that the value of the minimal $x_1$-term is $v(\binom zj X^j)$,
so it suffices to prove the {\it claim}: $v(\binom zj X^j)>v_0+2$ for $j\geq 4$.

Since $v_0=v(z(z-1))$, either $v_0=v(z)$ or $v_0=v(z-1)$.
Assume first that $v_0=v(z)$.
Using Kummer's formula it is easy to show that 
$v(\binom zj)=v(\binom{3^{v_0}}j)$ for $j\leq 3^{v_0}$.
Therefore if $j\leq 3^{v_0}$, $v(\binom zj)=v_0-v(j)$,
hence $v(\binom zjX^j)=v_0+j-v(j)$.
It is easy to see that $j-v(j)\geq 4$ when $j\geq 4$, hence
$v(\binom zj X^j)\geq v_0+4$ in this case, and we are done.
If $j>3^{v_0}$ then already $j>v_0+2$, so $v(\binom zj X^j)\geq j>v_0+2$,
and we are done.
This proves the claim for $v_0=v(z)$.

If $v_0=v(z-1)$, then $\binom zj=\frac z{z-j}\binom{z-1}j$
and $v(z)=0$, hence $v(\binom zj)=v(\binom{z-1}j)-v(z-j)$.
In particular, $v(\binom zj)=v(\binom{z-1}j)$ if and only if 
$j\equiv 0$ or $2(\mod 3)$.
But we have already shown the minimal $x_1$ term in this case has
value exceeding $v_0+2$ for $j\geq 4$, so we are done if $j\equiv 0$ or $2(\mod 3)$.
If $j\equiv 1(\mod 3)$ then since $v(\frac{z-j+1}j)=0$, 
$v(\binom zj)=v(\binom z{j-1})$,
so $v(\binom zj X^j)>v(\binom z{j-1}X^{j-1})$.
If $j>4$ then
since $j-1\equiv 0$ or $2(\mod 3)$, we are done by the previous case.
If $j=4$ then since $v(\binom z3 X^3)=v_0+2$,
we have the claim directly.
This finishes the claim.

Now suppose $k\geq 2$.
Set $Y=X-x_k 3^k$, then $v(Y)=1$.
By the Binomial Theorem,
$$
\tbinom zj X^j=\tbinom zj\sum_{i=0}^j\tbinom ji Y^{j-i} x_k^i 3^{ki}
$$
Evidently $x_k$ only appears in the $i\geq 1$ terms.
By inspection,
the minimal $x_k$-term of $\binom z2 X^2=\binom z2(Y+x_k 3^k)^2$ 
is $\binom z2 2x_1 x_k 3^{k+1}$.
To complete the proof it is enough to show 
that for $j\geq 3$ and $i\geq 1$,
$v(\tbinom zj)+v(\tbinom ji)+j-i+ki>v_0+k+1$, 
or 
$$
v(\tbinom zj)+v(\tbinom ji)+j+(i-1)(k-1)>v_0+2.
$$
If $j>3^{v_0}$ then since $v_0\geq 1$ we have $j>v_0+2$, 
and we are done.
Assume for the rest of the proof that $j\leq 3^{v_0}$. 
By Kummer's formula again, $v(\binom zj)=v_0-v(j)$ 
if $v_0=v(z)$ or if $v_0=v(z-1)$ and $j\equiv 0$ or $2(\mod 3)$,
and we are done in these cases if $j-v(j)>2$.
In particular, we are done if $j\geq 4$.
If $j=3\leq 3^{v_0}$ and $i\in\{1,2\}$ then $v(\binom ji)=1$ and we are done.
If $j=3\leq 3^{v_0}$ and $i=3$ 
then the left side of the inequality becomes $v_0-1+3+2(k-1)=v_0+2+2(k-1)$,
and since $k>1$, we are done in this case.
Finally, suppose $v_0=v(z-1)$ and $j\equiv 1(\mod 3)$.
Then $j\geq 4$.
As before, $v(\frac{z-j+1}j)=0$ implies $v(\binom zj)=v(\binom z{j-1})$.
Since $j-1\leq 3^{v_0}$ and $j-1\equiv 0(\mod 3)$, Kummer's formula yields
$v(\binom zj)=v_0-v(j-1)$, and we are done if $j-v(j-1)>2$.
Since $j\geq 4$, this completes the proof.

\endpf

The next result proves the existence of $q$ such that
$[\i_q]_{3^n}$ fixes a given $z\in\Z_3$.

\proclaim{Proposition 5.2 ``Existence of $q$''}
Suppose $n\geq 3$, $z\in\Z_3$, $v_0=v(z(z-1))$, and $1\leq v_0\leq n-2$.
Then there exists an element 
$q_{n-v_0}=1+a_1 3+\cdots +a_{n-v_0-1}3^{n-v_0-1}\in U^{(1)}-U^{(2)}$ such that
$\i_{q_{n-v_0}}(z)\equiv z(\mod 3^n)$.
If $q\in U^{(1)}$, then
$\i_q(z)\equiv z(\mod 3^n)$ if and only if $q\equiv q_{n-v_0}(\mod 3^{n-v_0})$.
If $v(z)\geq 1$ then $q_{n-v_0}\equiv 7(\mod 9)$, and if $v(z-1)\geq 1$ then 
$q_{n-v_0}\equiv 4(\mod 9)$.
\endproclaim

\Pf
Let $X=x_1 3+x_2 3^2+\cdots\in 3\Z[x_1,x_2,\dots]_3$.
We will show that the equation $S=\sum_{j\geq 2}\binom zj X^j\equiv 0(\mod 3^{n+1})$
has a solution $X=A\in\Z_3$ of value $1$,
that is uniquely determined modulo $3^{n-v_0}$.
Equivalently, $(1+X)^z\equiv 1+zX(\mod 3^{n+1})$ has such a solution, and setting $q=A+1$,
we obtain an element $q\in U^{(1)}-U^{(2)}$
satisfying $\i_q(z)\equiv z(\mod 3^n)$, uniquely determined $(\mod 3^{n-v_0})$.
We then set $q_{n-v_0}=1+a_1 3+\cdots+a_{n-v_0-1}3^{n-v_0-1}$.
We will treat the $v_0=v(z)$ case in detail, 
then indicate the changes needed to prove the $v_0=v(z-1)$ case.

We proceed inductively, showing that for all $n=m+v_0+1$,
$\sum_{j\geq 2}\binom zj X^j\equiv 0(\mod 3^{n+1})$
uniquely determines $a_1,\dots,a_m$.
Suppose $m=1$, so $n=v_0+2$.
By Lemma 5.1 the minimal $x_1$-term of $S$ is
$$
\tbinom z2 x_1^2 3^2+\tbinom z3 x_1^3 3^3=
\tfrac 12(z-1)x_1^2(1+(z-2)x_1)3^{v_0+2}.
$$ 
Making this $0(\mod 3^{v_0+3})$
is the same as solving
$\frac 12(z-1)x_1^2(1+(z-2)x_1)\equiv 0(\mod 3)$.
Since $z-1$ is invertible $(\mod 3)$ and $z\equiv 0(\mod 3)$,
the solutions are $x_1=0$ and $x_1=2$.
We eliminate $x_1=0$, since we require $q\in U^{(1)}-U^{(2)}$,
and we are left with a unique solution $x_1=a_1=2$.
This completes the base case.

Suppose we have uniquely determined $x_i=a_i$ in $S$ for $i=1,\dots,m$, 
so that $S\equiv 0(\mod 3^{n+1})$, where $n=m+v_0+1$. 
We will show $S\equiv 0(\mod 3^{n+2})$ uniquely determines $x_{m+1}=a_{m+1}$.
By Lemma 5.1 the minimal $x_{m+1}$-term of $S$ is 
$\binom z2 2x_1 x_{m+1}3^{m+2}=(z-1)x_1 x_{m+1}3^{v_0+m+2}$.
Thus $S(\mod 3^{v_0+m+3})$ is linear in $x_{m+1}$,
and reduces to $f(a_1,\dots,a_m)+(z-1)a_1 x_{m+1}\equiv 0(\mod 3)$
for some polynomial $f\in\Z[x_1,\dots,x_m]$.
Since $z-1$ and $a_1=2$
are invertible $(\mod 3)$, there is a unique solution $x_{m+1}=a_{m+1}\in\{0,1,2\}$,
as desired.
By induction, for all $n$,
the equation $S\equiv 0(\mod 3^{n+1})$ determines the coefficients $a_1,\dots,a_{n-v_0-1}$. 

The proof shows that if
$q_{n-v_0}=1+2\cdot 3+a_2 3^2+\cdots+a_{n-v_0-1}3^{n-v_0-1}$,
then for all $q\in U^{(1)}-U^{(2)}$, $\i_q(z)\equiv z(\mod 3^n)$ if and only if
$q\equiv q_{n-v_0}(\mod 3^{n-v_0})$.
This completes the proof of the $v_0=v(z)$ case.

If $v_0=v(z-1)$ then
by Lemma 5.1 the minimal $x_k$-terms have the same form with minor modifications, 
and the same argument uniquely determines the coefficients $a_1,\dots,a_{n-v_0-1}$.
However, $x_1=a_1$ is determined by the equation 
$\frac 12 zx_1^2(1+(z-2)x_1)\equiv 0(\mod 3)$, and since $z\equiv 1(\mod 3)$ this
forces $a_1=1$.
This completes the proof of the $v_0=v(z-1)$ case.

\endpf

As a corollary
we obtain a class of fixed points that do not appear in pairs,
unlike the others. 

\proclaim
{Corollary 5.3}
Fix a number $n\geq 3$.
\roster
\item"{a.}"
If $q\equiv 7(\mod 9)$ then $c\cdot 3^{n-2}$ is a fixed point for $[\i_q]_{3^n}$
for all $c\in\Z_3$.
\item"{b.}"
If $q\equiv 4(\mod 9)$ then $c\cdot 3^{n-2}+1$ is a fixed point for $[\i_q]_{3^n}$
for all $c\in\Z_3$.
\endroster
\endproclaim

\Pf
If $v(z(z-1))>n-2$ then $z$ is fixed by Theorem 4.2.
If $v(z(z-1))=n-2$, then by Proposition 5.2, 
for all $q\in U^{(1)}-U^{(2)}$,
$[\i_q]_{3^n}$ fixes $z$ if and only if $q\equiv 1+a_1 3(\mod 3^2)$.
Since $a_1=2$ when $q\equiv 7(\mod 9)$ and $a_1=1$ when $q\equiv 4(\mod 9)$, this
completes the proof.

\endpf

We now pass to the limit to assign a 
unique $q\in U^{(1)}-U^{(2)}$
to every element 
$z\in\Z_3$ such that $1\leq v(z(z-1))<\infty$.

\proclaim
{Corollary 5.4}
Suppose $z\in 3\Z_3\cup(1+3\Z_3)-\{0,1\}$. 
Then there exists a unique $q\neq 1$ in $U^{(1)}$ such that
$z$ is a $3$-adic fixed point for $\i_q$.
If $z\equiv 0(\mod 3)$ then $q\equiv 7(\mod 9)$, and if $z\equiv 1(\mod 3)$ then $q\equiv 4(\mod 9)$.
\endproclaim

\Pf
Let $v_0=v(z(z-1))$.
Since $z\neq 0,1$, we have $1\leq v_0<\infty$, 
and we may apply Proposition 5.2, for $n\geq v_0+2$.
Therefore for all $n\geq v_0+2$,
let $q_{n-v_0}$ denote the number determined for $z$ modulo $3^n$
by Proposition 5.2.
The $q_{n-v_0}$ are compatible,
in the sense that for all $m$ and $n$ such that $m<n$, 
$q_{n-v_0}\equiv q_{m-v_0}(\mod 3^{m-v_0})$.
For since $\i_{q_{n-v_0}}(z)\equiv z(\mod 3^n)$, we have 
$\i_{q_{n-v_0}}(z)\equiv z(\mod 3^m)$ for all $m<n$, therefore by Proposition 5.2,
$q_{n-v_0}\equiv q_{m-v_0}(\mod 3^{m-v_0})$, as desired.
Thus the number $q=\lim_{n\to\infty}q_{n-v_0}$ is well defined,
and since $q\equiv q_{n-v_0}(\mod 3^{n-v_0})$ for all $n$, $\i_q(z)\equiv z(\mod 3^n)$
for all $n$, hence $z$ is a $3$-adic fixed point for $\i_q$.

\endpf

\proclaim
{Lemma 5.5}
Fix a number $n\geq 3$.
Suppose $q\in U^{(1)}-U^{(2)}$, and $a_0\in 3\Z_3$, $a_0\neq 0$.
\roster
\item"{a.}"
If $a_0$ is a fixed point for $[\i_q]_{3^n}$
then $ca_0$ is a fixed point for $[\i_q]_{3^n}$ for $c\in\Z_3$
if and only if $c(c-1)\in 3^{n-2v(a_0)-1}\Z_3$.
\item"{b.}"
If $a_0+1$ is a fixed point for $[\i_q]_{3^n}$
then $ca_0+1$ is a fixed point for $[\i_q]_{3^n}$ for $c\in\Z_3$
if and only if $c(c-1)\in 3^{n-2v(a_0)-1}\Z_3$.
\endroster
\endproclaim

\Pf
Let $v_0=v(a_0)$.
Suppose $a_0$ is a fixed point.
By Lemma 1.2, $\i_q(ca_0)=\i_{q^{a_0}}(c)\i_q(a_0)$,
so $ca_0$ is a fixed point if and only if
$\i_{q^{a_0}}(c)\i_q(a_0)\equiv ca_0(\mod 3^n)$. 
Since $v(a_0)=v_0$ and $\i_q(a_0)\equiv a_0(\mod 3^n)$,
we see that $ca_0$ is a fixed point if and only if
$\i_{q^{a_0}}(c)\equiv c(\mod 3^{n-v_0})$.
To prove (b), suppose $a_0+1$ is a fixed point.
By the cocycle condition, $\i_q(ca_0+1)\equiv ca_0+1(\mod 3^n)$ is equivalent to
$q\i_q(ca_0)\equiv ca_0(\mod 3^n)$, hence
$ca_0+1$ is a fixed point if and only if
$q\i_{q^{a_0}}(c)\i_q(a_0)\equiv ca_0(\mod 3^n)$.
Since $\i_q(a_0+1)\equiv a_0+1(\mod 3^n)$ we have $q\i_q(a_0)\equiv a_0(\mod 3^n)$,
hence $ca_0+1$ is a fixed point if and only if 
$\i_{q^{a_0}}(c)\equiv c(\mod 3^{n-v_0})$.

Since $q\in U^{(1)}$ and $v_0\neq 0$, $q^{a_0}\in U^{(2)}$, 
and Theorem 4.2 shows $c$ is a fixed
point for $[\i_{q^{a_0}}]_{3^{n-v_0}}$ if and only if $c$ or $c-1$ is a multiple
of $o([q^{a_0}]_{3^{n-v_0}})=o([q]_{3^{n-2v_0}})$.
But since $q\in U^{(1)}-U^{(2)}$, $o([q]_{3^{n-2v_0}})=3^{n-2v_0-1}$ by Lemma 2.2, 
as claimed.

\endpf

\proclaim
{Definition 5.6}
{\rm
Fix a number $n\geq 2$.
Suppose $q\in U^{(1)}-U^{(2)}$ and 
$z_0$ is a modular fixed point of $[\i_q]_{3^n}$ such that
$v_0:=v(z_0(z_0-1))$ has the smallest nonzero value.
The {\it period} of the fixed point set for $[\i_q]_{3^n}$
is the number $\tau=3^{n-v_0-1}$.
A fixed point $z$ for $[\i_q]_{3^n}$ is called
a {\it rooted} fixed point if $1\leq v(z(z-1))<\frac{n-1}2$, and 
a {\it drifting} fixed point if $v(z(z-1))\geq\frac{n-1}2$.
}
\endproclaim

If $z\not\in\{0,1\}$ is a $3$-adic fixed point for $\i_q$, 
then for $n>2v(z(z-1))+1$, $z$ is a rooted fixed point for $[\i_q]_{3^n}$.
We will show below that conversely a rooted fixed point indicates the existence
of a $3$-adic fixed point.

\flushpar
{\bf Summary 5.7.}
We summarize the situation so far.
By Theorem 4.2, the subgroup $o([q]_{3^n})\Z_3$ along with its coset
$1+o([q]_{3^n})\Z_3$ form a subset of fixed points for $[\i_q]_{3^n}$,
the ``modular fixed point pairs''.
This set can be computed directly from $q$.
In fact, by Lemma 2.2, $o([q]_{3^n})=3^{n-m_0}$, where $m_0=v_3(q-1)$.
If $q\in U^{(2)}$, then these are the only fixed points.
If $q\equiv 4$ or $7(\mod 9)$, then the set $1+3^{n-2}\Z_3$ or $3^{n-2}\Z_3$,
respectively, give additional fixed points, by Corollary 5.3.
We will now see that the remaining fixed points are much more obscure.

The next result
gives the complete modular fixed point set structure,
{\it given} a fixed point $z_0$ for $[\i_q]_{3^n}$
such that $v_0:=v(z_0(z_0-1))$ is minimal.
We find two distinct cases, depending on whether this $v_0$
is less than $\frac{n-1}2$, i.e., whether there exists a rooted fixed point.
If $z_0$ is a rooted fixed point,
then every $z$ congruent to $z_0$ modulo $\tau=3^{n-v_0-1}$ is also a rooted fixed point,
and this is the complete rooted fixed point set.
This set is irregular in the sense that the
valuation data of a given $z$ does not by itself predict membership;
the number has to have a certain residue.
Aside from rooted fixed points, there is a set of drifting fixed points,
determined by valuation data alone.

\proclaim
{Theorem 5.8}
Fix a number $n\geq 2$.
Suppose $q\in U^{(1)}-U^{(2)}$. 
Let $z_0\in\Z_3$ be a modular fixed point of $[\i_q]_{3^n}$ such that
$v_0:=v(z_0(z_0-1))$ has the smallest nonzero value.
Note that by Theorem 4.2, $v_0\leq n-1$.
Let $\tau=3^{n-v_0-1}$.
\roster
\item"{a.}"
Suppose $q\equiv 7(\mod 9)$.
The fixed point set for $[\i_q]_{3^n}$ is
$$
\cases
(z_0+\tau\Z_3)\cup\tau\Z_3\cup(1+3^{n-1}\Z_3)
&\text{if $v_0<\frac{n-1}2$}\\
3^{\lfloor\frac n2\rfloor}\Z_3\cup(1+3^{n-1}\Z_3)
&\text{if $v_0\geq\frac{n-1}2$}\\
\endcases
$$
\item"{b.}"
Suppose $q\equiv 4(\mod 9)$.
The fixed point set for $[\i_q]_{3^n}$ is
$$
\cases
(z_0+\tau\Z_3)\cup(1+\tau\Z_3)\cup 3^{n-1}\Z_3
&\text{if $v_0<\frac{n-1}2$}\\
(1+3^{\lfloor\frac n2\rfloor}\Z_3)\cup 3^{n-1}\Z_3
&\text{if $v_0\geq\frac{n-1}2$}\\
\endcases
$$
\endroster
\endproclaim

\Pf
Set $o_{3^n}=o([q]_{3^n})$.
If $n=2$ then $o_{3^n}=3$, and by Theorem 4.2 the fixed point set is
$3\Z_3\cup (1+3\Z_3)$ for any $q\in U^{(1)}-U^{(2)}$, $q\neq 1$.
Since here we have $v_0=1\geq\frac{n-1}2$, we obtain the desired expression.

Suppose $n\geq 3$ and $q\equiv 7(\mod 9)$.
Then $v_0=v(z_0)$.
By Corollary 5.3, $1\leq v_0\leq n-2$.
By Proposition 5.2, any fixed point for $\i_q$ with $1\leq v_0\leq n-2$
is divisible by $3$.
We first treat the case $v_0<\frac{n-1}2$.
Since $z_0\in 3\Z_3$,
$cz_0$ is a fixed point if and only if $c(c-1)\in 3^{n-2v_0-1}\Z_3$,
by Lemma 5.5(a).
If $c-1\in 3^{n-2v_0-1}\Z_3$, then 
$cz_0\in z_0+z_03^{n-2v_0-1}\Z_3=z_0+3^{n-v_0-1}\Z_3=z_0+\tau\Z_3$.
If $c\in 3^{n-2v_0-1}\Z_3$ then similarly $cz_0\in \tau\Z_3$.
Thus the fixed point set contains $(z_0+\tau\Z_3)\cup\tau\Z_3$.
Any remaining fixed points $z:v(z)\leq n-2$
cannot be multiples of $z_0$ by Lemma 5.5,
hence by the minimality of the value of $z_0$ they must have $3$-value $0$.  
But by Proposition 5.2, all numbers $z: 1\leq v(z-1)\leq n-2$
can only be fixed points for $q\equiv 4(\mod 9)$,
and by Theorem 4.2, all numbers $z\equiv 2(\mod 3)$ are not fixed points 
for any $q\neq 1$.
On the other hand, by Theorem 4.2
every element in $1+3^{n-1}\Z_3$ is a fixed point for $[\i_q]_{3^n}$.
This gives the desired fixed point set for the case $v_0<\frac{n-1}2$
and $q\equiv 7(\mod 9)$.
The argument for $v_0<\frac{n-1}2$ and $q\equiv 4(\mod 9)$ is similar, 
and we just indicate the changes.
Let $a_0=z_0-1$, so $v(a_0)=v_0$.
By Lemma 5.5(b), if $v_0<\frac{n-1}2$ then $ca_0+1$ is fixed if and only
if $c(c-1)\in 3^{n-1-2v_0}\Z_3$, and
we see the fixed point set contains $(z_0+\tau\Z_3)\cup (1+\tau\Z_3)$.
The only remaining fixed points are given by the set $3^{n-1}\Z_3$ by Theorem 4.2,
using Proposition 5.2 to rule out all $z:1\leq v(z)\leq n-2$, and Theorem 4.2
to rule out all $z:z\equiv 2(\mod 3)$.

Suppose $v_0\geq\frac{n-1}2$ and $q\equiv 7(\mod 9)$.
Then $2v_0\geq n-1$, so by Lemma 5.5(a),
$cz_0$ is a fixed point for all $c\in\Z_3$.
Any remaining fixed points cannot be multiples of $z_0$, and as above we add only
the set $1+3^{n-1}\Z_3$.
It remains to show that we can take $z_0=3^{\lfloor\frac n2\rfloor}$.
For this we need a lemma, which we will also use later.

\proclaim
{Lemma 5.9}
Fix numbers $n\geq 3$ and $v_0:1\leq v_0<\frac{n-1}2$.
The total number of generators $q\in U^{(1)}/U^{(n)}$
such that $[\i_q]_{3^n}$ has a rooted fixed point $z$ with $v(z(z-1))=v_0$ 
is $4\cdot 3^{n-v_0-2}$,
divided evenly between the different possibilities for $z(\mod 3^{n-v_0-1})$.
The total number of $q$
such that $[\i_q]_{3^n}$ has a rooted fixed point is 
$2(3^{n-2}-3^{\lfloor\frac{n-1}2\rfloor})$.
The number such that $[\i_q]_{3^n}$ has no rooted fixed point is
$2\cdot 3^{\lfloor\frac{n-1}2\rfloor}$.
In each case half of the number is for $q\equiv 7(\mod 9)$, the other half for
$q\equiv 4(\mod 9)$.
\endproclaim

\Pf
Suppose $q$ generates $U^{(1)}/U^{(n)}$.
By Definition 5.6,
$z\in\Z_3$ is a rooted fixed point for $[\i_q]_{3^n}$
if and only if $1\leq v(z(z-1))<\frac{n-1}2$.
By Theorem 5.8 so far, $[\i_q]_{3^n}$ has a rooted fixed point $z$
such that $v(z(z-1))=v_0$ if and only if it has a
(uniquely determined) rooted fixed point $z_0$ such that
$1<z_0<3^{n-v_0-1}$.
The number of possible distinguished rooted fixed points $z_0$
with $v(z_0)=v_0$
is thus the number of generators of $3^{v_0}\Z/3^{n-v_0-1}\Z$,
or $2\cdot 3^{n-2-2v_0}$.
The number with $v(z_0-1)=v_0$ is obviously the same,
so the total number of $z_0$ with $v(z_0(z_0-1))=v_0$ and $1<z_0<3^{n-v_0-1}$
is $4\cdot 3^{n-2-2v_0}$.
Conversely, by Proposition 5.2, each such $z_0$ defines a number
$q_{n-v_0}=1+a_1 3+\cdots+a_{n-v_0-1}3^{n-v_0-1}$, and 
among all generators $q\in U^{(1)}/U^{(n)}$,
$[\i_q]_{3^n}$ fixes $z_0$ if and only if $q\equiv q_{n-v_0}(\mod 3^{n-v_0})$.
Thus to each $z_0$ there are $3^{v_0}$ generators $q\in U^{(1)}/U^{(n)}$
such that $[\i_q]_{3^n}$ fixes $z_0$.
Therefore the total number of $q$ such that $[\i_q]_{3^n}$ 
has a rooted fixed point $z$ with
$v(z(z-1))=v_0$ is $4\cdot 3^{n-2-2v_0}\cdot 3^{v_0}=4\cdot 3^{n-v_0-2}$,
divided evenly between all of the possible $z_0(\mod 3^{n-v_0-1})$,
as claimed in the first statement.

It follows that for a given $n$, the total number of generators
$q\in U^{(1)}/U^{(n)}$ such that $[\i_q]_{3^n}$ has
a rooted fixed point is
$\sum_{v_0=1}^{\lceil\frac{n-1}2\rceil-1}4\cdot 3^{n-v_0-2}$,
and half of the $q$ are $4(\mod 9)$, the other half $7(\mod 9)$.
Now compute
$$
\align
\sum_{v_0=1}^{\lceil\frac{n-1}2\rceil-1}4\cdot 3^{n-v_0-2}
&=4\sum_{j=\lfloor\frac{n-1}2\rfloor}^{n-3} 3^j
=4\cdot 3^{\lfloor\frac{n-1}2\rfloor}
\sum_{j=0}^{n-\lfloor\frac{n-1}2\rfloor-3}3^j
\\ &
=4\cdot 3^{\lfloor\frac{n-1}2\rfloor}
\cdot\tfrac{3^{n-\lfloor\frac{n-1}2\rfloor-2}-1}{3-1}
\\ &
=2\cdot 3^{\lfloor\frac{n-1}2\rfloor}
\cdot(3^{n-\lfloor\frac{n-1}2\rfloor-2}-1)
\\ &
=2(3^{n-2}-3^{\lfloor\frac{n-1}2\rfloor}),
\endalign
$$
which proves the second statement.
The total number of generators $q\in U^{(1)}/U^{(n)}$
is $2\cdot 3^{n-2}$, since $a_1$ can assume one of the
two values $1,2$, and $a_2,\dots,a_{n-1}$ can assume any of the values $0,1,2$.
Therefore the number of such $q$ such that $[\i_q]_{3^n}$
has no rooted fixed points is
$2\cdot 3^{n-2}-2(3^{n-2}-3^{\lfloor\frac{n-2}2\rfloor})
=2\cdot 3^{\lfloor\frac{n-1}2\rfloor}$,
and again half are congruent to $7(\mod 9)$, the other half $4(\mod 9)$.
This completes the proof of the lemma.

\endpf

We continue with the proof of Theorem 5.8.
We are showing that if $n\geq 3$, $v_0\geq\frac{n-1}2$, and $q\equiv 7(\mod 9)$,
then $[\i_q]_{3^n}$ fixes $3^{\lfloor\frac n2\rfloor}$.

Suppose $n$ is odd.  Using Proposition 5.2, we count
the number of $q$ between $1$ and $3^n$
such that $[\i_q]_{3^n}$ fixes
$3^{\lfloor\frac n2\rfloor}=3^{\frac{n-1}2}$ and obtain $3^{\frac{n-1}2}$.
None of the $q$ for which $\i_q$ has a rooted fixed point can fix
$3^{\frac{n-1}2}$, since the rooted fixed point has value $v_0<\frac{n-1}2$,
and the corresponding period $\tau$ has value $n-v_0-1>\frac{n-1}2$,
leaving nothing in between.
Thus for all of the $3^{\frac{n-1}2}$ numbers $q$ between $1$ and $3^n$ 
for which $3^{\frac{n-1}2}$ is a fixed point for $[\i_q]_{3^n}$,
$\i_q$ has no rooted fixed points.
By Lemma 5.9, this accounts for all of the $\i_q$, $q\equiv 7(\mod 9)$,
with no rooted fixed points, as desired.

Suppose $n$ is even. 
We want to show that if $[\i_q]_{3^n}$ has no rooted fixed points 
then it fixes $3^{\frac n2}$.
To do it, we count the total number of $q$ for which $[\i_q]_{3^n}$
fixes $3^{\frac n2}$, and subtract the number that fix $3^{\frac n2}$
{\it and} have a corresponding rooted fixed point. 
We need to show that the result is the same as the number of $q$ for
which $[\i_q]_{3^n}$ has no rooted fixed points,
which by Lemma 5.9 is $3^{\lfloor\frac{n-1}2\rfloor}=3^{\frac n2-1}$.
By Theorem 5.8 so far, 
if $[\i_q]_{3^n}$ fixes $3^{\frac n2}$ and has a rooted fixed point,
we have $\frac n2\geq n-1-v_0$ and $v_0<\frac{n-1}2$, hence
this rooted fixed point must have value $v_0=\frac n2-1$.

Since $n+1$ is odd, by what we have just shown
there are $3^{\frac{n+1-1}2}=3^{\frac n2}$ 
numbers $q$ between $1$ and $3^{n+1}$ such that 
$3^{\lfloor\frac{n+1}2\rfloor}=3^{\frac n2}$ is a fixed point
for $[\i_q]_{3^{n+1}}$, and these $q$ account for all of
the $[\i_q]_{3^{n+1}}$ with no rooted fixed points.
Any fixed point for $[\i_q]_{3^{n+1}}$ is a fixed point for $[\i_q]_{3^n}$,
but the number of these $q$ that are distinct modulo $3^n$ is reduced by
a factor of $3$.
Thus we have $3^{\frac n2-1}=3^{\lfloor\frac{n-1}2\rfloor}$ 
numbers $q$ between $1$ and $3^n$ such that $[\i_q]_{3^n}$ fixes
$3^{\frac n2}$, {\it and} $[\i_q]_{3^n}$ is descended from 
$[\i_q]_{3^{n+1}}$ that have no rooted fixed points.

We {\it claim} that all of these $[\i_q]_{3^n}$ have no rooted fixed points.
To prove it, we will show that all of the $[\i_q]_{3^n}$ that
fix $3^{\frac n2}$ and {\it have} rooted fixed points 
are descended from $[\i_q]_{3^{n+1}}$ that have rooted fixed points.
By Lemma 5.9, there are $2\cdot 3^{\frac n2}$ numbers $q$ such that $[\i_q]_{3^{n+1}}$
has a rooted fixed point of value $v_0=\frac n2-1$.
These rooted fixed points for $[\i_q]_{3^{n+1}}$
are rooted fixed points for $[\i_q]_{3^n}$, of value 
$\frac n2-1<\frac{n-1}2$, so there are $2\cdot 3^{\frac n2-1}$ of these $q$ such that
$[\i_q]_{3^n}$ has a rooted fixed point of value $v_0=\frac n2-1$.
By Lemma 5.9, this is all of the $q$ between $1$ and $3^n$
such that $[\i_q]_{3^n}$ has a rooted fixed point of value $\frac n2-1$.
The period corresponding to these rooted fixed points is $3^{n-v_0-1}=3^{\frac n2}$,
and each of these $[\i_q]_{3^n}$ fixes the period.
This proves the claim.
Thus we have $3^{\frac n2-1}$ numbers $q$ between $1$ and $3^n$ for which
$[\i_q]_{3^n}$ fixes $3^{\frac n2}$, such that $[\i_q]_{3^n}$
has no rooted fixed points.
By Lemma 5.9, this accounts for all such $q$, as desired.

The argument for $v_0\geq\frac{n-1}2$ and $q\equiv 4(\mod 9)$ is similar.
Since $v_0\geq\frac{n-1}2$, by Lemma 5.5(b) we obtain
the set $z_0\Z_3\cup 3^{n-1}\Z_3$, as desired;
we then apply Lemma 5.9 as before to show that we can take 
$z_0=1+3^{\lfloor\frac n2\rfloor}$.
This completes the proof.

\endpf

We now count the number of fixed points in $\Z_p(\mod p^n)$
for $[\i_q]_{p^n}$, where $p$ is any prime and $q\in U^{(1)}$.

\proclaim
{Corollary 5.10 ``Fixed Point Count''}
Suppose $q\in U^{(1)}$, $n\in\N$, and $m_0=v_p(q-1)$.
If $p=3$ and $q\in U^{(1)}-U^{(2)}$, let $v_0$ be the smallest nonzero value
of $v(z(z-1))$ for any modular fixed point $q$ of $[\i_q]_{p^n}$.
Then
if $q\equiv 1(\mod p^n)$ then every point of $\Z_p$ is a fixed point of $[\i_q]_{p^n}$.
If $q\not\equiv  1(\mod p^n)$
then the number of fixed points $z$ modulo $p^n\Z_p$ of $[\i_q]_{p^n}$ is
$$
\cases
2p^{m_0} &\text{if $p\neq 3$ or $q\in U^{(2)}$, and $p$ is odd}\\
2^{m_0} &\text{if $p=2$}\\
2\cdot 3^{v_0+1}+3 &\text{if $p=3$, $q\in U^{(1)}-U^{(2)}$, and $v_0<\tfrac{n-1}2$}\\
3^{n-\lfloor\frac n2\rfloor}+3 &\text{if $p=3$, $q\in U^{(1)}-U^{(2)}$,
and $v_0\geq\tfrac{n-1}2$}\\
\endcases
$$
\endproclaim

\Pf
The first statement is clear.
If $q\not\equiv  1(\mod p^n)$ and $q\in U^{(2)}$,
then $o_{p^n}=p^{n-m_0}$ by Lemma 2.2.
The rest is a simple count, using Theorem 4.2 if $p\neq 3$ or $q\in U^{(2)}$, 
and Theorem 5.8 if $p=3$ and $q\in U^{(1)}-U^{(2)}$.  
In the first instance the fixed point set modulo $p^n\Z_p$
has the form $a_0\Z/p^n\Z\cup(1+a_0\Z/p^n\Z)$,
where $a_0=o_{p^n}=p^{n-m_0}$ if $p$ is odd, $a_0=2\cdot o_{2^n}$ 
if $p=2$ and $q\in U^{(2)}$,
and $a_0=2^n$ if $p=2$ and $q\in U^{(1)}-U^{(2)}$.
Thus the count is $2\cdot p^n/(p^{n-m_0})=2\cdot p^{m_0}$, 
$2\cdot 2^n/(2\cdot 2^{n-m_0})=2^{m_0}$, and $2$, respectively.
In the second instance, with $p=3$, we have 
$(2\cdot 3^n/\tau)+3=2\cdot 3^{n-(n-v_0-1)}+3=2\cdot 3^{v_0+1}+3$
if $q\in U^{(1)}-U^{(2)}$ and $v_0<\frac{n-1}2$, 
and $(3^n/3^{\lfloor\frac n2\rfloor})+3=3^{n-\lfloor\frac n2\rfloor}+3$ 
if $q\in U^{(1)}-U^{(2)}$ and $v_0\geq\frac{n-1}2$.

\endpf

\flushpar
{\bf Remark 5.11.}
We will show in Corollary 6.4 that the number of fixed points for $[\i_q]_{p^n}$
is asymptotically stable for all $q\in U^{(1)}$ and primes $p$, with the exception
of exactly two values of $q$ for the prime $p=3$.

We will now show how to use the 
rooted fixed points to construct $3$-adic fixed points.

\proclaim
{Lemma 5.12 ``Propagation of Rooted Fixed Points''}
Suppose $n\geq 3$, $q\in U^{(1)}-U^{(2)}$,
$z_0\in\Z_3$, and $v_0:=v(z_0(z_0-1))<\frac{n-1}2$.
If $\i_q(z_0)\equiv z_0(\mod 3^n)$ 
then there exists a unique $c\in\{0,1,2\}$ such that 
$\i_q(z_0+c 3^{n-v_0-1})\equiv z_0+c 3^{n-v_0-1}(\mod 3^{n+1})$.
In particular, if $[\i_q]_{3^{n_0}}$ has a rooted fixed point for some $n_0$,
then $[\i_q]_{3^n}$ has a rooted fixed point for all $n\geq n_0$.
\endproclaim

\Pf
Suppose $\i_q(z_0)\equiv z_0(\mod 3^n)$.
If $[\i_q]_{3^{n+1}}$ has a rooted fixed point $z_0'$ 
satisfying $v(z_0'(z_0'-1))<\frac{n-1}2$,
then by Definition 5.6, $z_0'$ is a rooted fixed point for $[\i_q]_{3^n}$.  
By the uniqueness of $v_0$ in Theorem 5.8, we have
$v(z_0'(z_0'-1))=v_0$, and by Theorem 5.8, $z_0'=z_0+c3^{n-v_0-1}$
for some $c\in\{0,1,2\}$, as desired.  Since by Theorem 5.8, $z_0'$
is uniquely determined modulo $3^{(n+1)-1-v_0}=3^{n-v_0}$, $c$ is
uniquely determined.

We may assume for the rest of the proof that
$[\i_q]_{3^{n+1}}$ either has no rooted fixed points, or it has a rooted fixed
point $z_0'$ satisfying $v_0':=v(z_0'(z_0'-1))\geq\frac{n-1}2$.
Thus if $z_0'$ is a rooted fixed point for $[\i_q]_{3^{n+1}}$, by definition 
$v_0'<\frac{(n+1)-1}2=\frac n2$, so $v_0'=\frac{n-1}2$ is forced.
Let $\tau=3^{n-v_0-1}$.
We {\it claim} that $c\tau$ or $c\tau+1$ is a fixed point for $[\i_q]_{3^{n+1}}$ for any
$c\in\Z_3$, depending, as usual, on whether $q\equiv 7$ or $4(\mod 9)$.
We will only prove it for $q\equiv 7(\mod 9)$; the $q\equiv 4(\mod 9)$ case is similar.
If $[\i_q]_{3^{n+1}}$ has no rooted fixed point then by Theorem 5.8, 
$[\i_q]_{3^{n+1}}$ fixes every multiple of 
$3^{\lfloor\frac{n+1}2\rfloor}$.
Since 
$\tau=3^{n-v_0-1}$ and $v_0<\frac{n-1}2$, we have $v(\tau)>\frac{n-1}2$, hence
$v(\tau)\geq 3^{\lfloor\frac{n+1}2\rfloor}$.  Therefore
$c\tau$ is fixed by $[\i_q]_{3^{n+1}}$, for every $c\in\Z_3$, as desired.
If $[\i_q]_{3^{n+1}}$ has a rooted fixed point $z_0'$, so
$v_0'=\frac{n-1}2$, then
$n$ is odd, and by Theorem 5.8, $[\i_q]_{3^{n+1}}$ fixes every multiple of 
$\tau'=3^{(n+1)-1-v_0'}=3^{\frac{n+1}2}$.
Since $v_0<\frac{n-1}2$ and $n$ is odd, $\tau=3^{n-v_0-1}\geq 3^{\frac{n+1}2}$,
so $\tau$ is a multiple of $\tau'$.
Therefore $[\i_q]_{3^{n+1}}$ fixes every multiple of $\tau$.  
This proves the claim.

Since $\i_q(z_0)\equiv z_0(\mod 3^n)$, we have $\i_q(z_0)\equiv z_0+a_n 3^n(\mod 3^{n+1})$
for some $a_n\in\{0,1,2\}$.  Let $X=q-1$, then $v(X)=1$.
If $q\equiv 7(\mod 9)$ then $v(z_0\tau X)=n$, hence $z_0\tau X\equiv b_n 3^n(\mod 3^{n+1})$
for some $b_n\in\{1,2\}$; if $q\equiv 4(\mod 9)$ then similarly we have 
$(z_0-1)\tau X\equiv b_n 3^n(\mod 3^{n+1})$.  In either case, set $c=-a_n b_n$.
Note since $b_n$ is invertible modulo $3$, we have $cb_n\equiv -a_n(\mod 3)$.
We {\it claim} $\i_q(z_0+c\tau)\equiv z_0+c\tau(\mod 3^{n+1})$.

Assume first that $q\equiv 7(\mod 9)$, so $v_0=v(z_0)$.
By the Binomial theorem,
$$
q^{z_0}=1+z_0 X+\sum_{j=2}^\infty\tbinom{z_0}j X^j.
$$
By Lemma 5.1, $v(\sum_{j=2}^\infty\tbinom{z_0}j X^j)\geq v_0+2$,
consequently $q^{z_0}c\tau\equiv (1+z_0 X)c\tau(\mod 3^{n+1})$.
By the previous claim, $\i_q(c\tau)\equiv c\tau(\mod 3^{n+1})$.
Therefore we compute modulo $3^{n+1}$,
$$
\align
\i_q(z_0+c\tau)
&\equiv \i_q(z_0)+q^{z_0}\i_q(c\tau)
\equiv z_0+a_n 3^n+q^{z_0}c\tau
\equiv z_0+a_n 3^n+(1+z_0 X)c\tau
\\ &
\equiv z_0 + c\tau+a_n 3^n + cb_n 3^n
\equiv z_0+c\tau\;(\mod 3^{n+1})
\endalign
$$
as desired.
If $q\equiv 4(\mod 9)$, then $v_0=v(z_0-1)$.
Since $\i_q(z_0)\equiv z_0+a_n 3^n(\mod 3^{n+1})$, we have
$q\i_q(z_0-1)=z_0-1+a_n 3^n$ by the cocycle condition,
and since $\i_q(c\tau+1)\equiv c\tau+1(\mod 3^{n+1})$,
we have $q\i_q(c\tau)\equiv c\tau(\mod 3^{n+1})$.
By the Binomial theorem, 
$q^{z_0-1}c\tau\equiv (1+(z_0-1)X)c\tau(\mod 3^{n+1})$, as before. 
Therefore we compute modulo $3^{n+1}$ using the cocycle condition,
$$
\align
\i_q(z_0+c\tau)
&\equiv 1+q\i_q(z_0-1+c\tau)
\equiv 1+q(\i_q(z_0-1)+q^{z_0-1}\i_q(c\tau))
\equiv z_0+a_n 3^n+q^{z_0-1}c\tau
\\ &
\equiv z_0+a_n 3^n+(1+(z_0-1)X)c\tau
\equiv z_0+a_n 3^n+c\tau+cb_n 3^n
\equiv z_0+c\tau\;(\mod 3^{n+1})
\endalign
$$
as desired.
This proves the claim.
It remains to show the solution $c\in\{0,1,2\}$ is unique.  But since $z_0+c\tau$
is a rooted fixed point for $[\i_q]_{3^{n+1}}$, its residue is uniquely determined modulo
$3^{(n+1)-1-v_0}=3^{n-v_0}$, and since $\tau=3^{n-v_0-1}$,
$c$ is uniquely determined as an element of $\{0,1,2\}$.

To prove the last statement, suppose $[\i_q]_{3^{n_0}}$ has a rooted fixed point.
By Definition 5.6, this is a fixed point $z\in\Z_3$ such that 
$v_0:=v(z(z-1))<\frac{n_0-1}2$.  By the above, $[\i_q]_{3^{n_0+1}}$ 
has a fixed point
with the same value $v_0$, and evidently $v_0<\frac{(n_0+1)-1}2=\frac{n_0}2$.
Therefore $z$ is a rooted fixed point for $[\i_q]_{3^{n_0+1}}$, and by
induction $[\i_q]_{3^n}$ has a rooted fixed point for all $n\geq n_0$.
This completes the proof.

\endpf

\proclaim
{Theorem 5.13 ``Existence of $3$-adic Fixed Points''}
Let $q\in U^{(1)}-U^{(2)}$.
Suppose $z$ is a rooted fixed point for $[\i_q]_{3^n}$ for some $n$,
and $v_0=v(z(z-1))$.
Then there exists a unique $3$-adic fixed point $z_q$ for $\i_q$,
and $z_q\equiv z(\mod 3^{n-v_0-1})$.
In particular, $v(z_q(z_q-1))=v_0$.
\endproclaim

\Pf
We will construct a $3$-adic fixed point $z_q=\lim_{i\to\infty}z_i$.
Let $z_{n-v_0-1}$ be the residue of $z(\mod 3^{n-v_0-1})$.
By Theorem 5.8, $z_{n-v_0-1}$ is a rooted fixed point of $[\i_q]_{3^n}$. 
Thus $z_{n-v_0-1}$ is the unique rooted fixed point for $[\i_q]_{3^n}$
between $1$ and the period $3^{n-v_0-1}$.
Set $z_1=\cdots=z_{n-v_0-1}$.
As for all rooted fixed points,
$v(z_{n-v_0-1}(z_{n-v_0-1}-1))=v_0$. 
By Lemma 5.12 there exists a uniquely defined number $c\in\{0,1,2\}$ such that
$[\i_q]_{3^{n+1}}$ fixes
$z_{n-v_0}:=z_{n-v_0-1}+c 3^{n-v_0-1}$.
Note $z_{n-v_0}\equiv z_{n-v_0-1}(\mod 3^{n-v_0-1})$, and $z_{n-v_0}$ is between
$1$ and $3^{n-v_0}$.
Since $v_0<n-v_0-1$, we have $v(z_{n-v_0}(z_{n-v_0}-1))=v_0$.
In this way we may 
define numbers $z_m$ inductively for all $m\in\N$, with the property that
$z_m$ is a fixed point for $[\i_q]_{3^{m+1+v_0}}$, 
$z_{m+1}\equiv z_m(\mod 3^m)$, and $v(z_m(z_m-1))=v_0$.
If $r>m$ then $z_r-z_m\in 3^m\Z_3$, hence
we have a well defined $3$-adic number $z_q=\lim_{m\to\infty}z_m\in\Z_3$,
such that $z_q\equiv z_m(\mod 3^m)$.
By construction, $z_q$ is a fixed point for $[\i_q]_{3^n}$ for all $n$, 
therefore $z_q$ is a $3$-adic fixed point for $\i_q$, and 
$z_q=z_{n-v_0-1}\equiv z(\mod 3^{n-v_0-1})$.
In particular, since $v_0<\frac{n-1}2<n-v_0-1$,
$v(z_q(z_q-1))=v(z(z-1))=v_0$.
This completes the proof.

\endpf

\flushpar
{\bf Remark 5.14}
By Theorem 5.8,
the rooted fixed points $z$ for $[\i_q]_{3^n}$ for a given $q\in U^{(1)}-U^{(2)}$
are distinguished from the set of drifting modular fixed points
by their lack of regularity:  they are the only ones that do not repeat in
intervals of $3^{v(z(z-1))}$.  
If, for example, $z=4$ is a rooted fixed point,
then since $v(4(4-1))=1$, $4+3=7$ is not.  
Thus if the fixed point set is observed to be regular, then there are no
rooted fixed points.
By Theorem 5.13 and its proof, we {\it need} 
the rooted fixed points to construct the $3$-adic fixed point, via residues. 
Even then, the rooted fixed points for $[\i_q]_{3^n}$ 
only give us the $3$-adic fixed point residues modulo $3^{n-v_0-1}$;
in general we cannot deduce them modulo $3^n$.
On the other hand, if we know $z_q$ is the $3$-adic fixed point for $q$
then we easily find all of the rooted fixed points for $[\i_q]_{3^n}$, 
for sufficiently large $n$, 
by taking the residue of $z_q$ modulo $3^n$, and adding to it all multiples
of the corresponding period $3^{n-v_0-1_q(z_q-1))}$.

We next aim to define a homeomorphism
from $3\Z_3\cup(1+3\Z_3)$ to $U^{(1)}-U^{(2)}$
sending $z$ to the unique $q$ such that $\i_q(z)=z$.
If $0<v(z(z-1))<\infty$, then the assignment is well defined by Theorem 5.13.
This leaves out $z=0$ and $z=1$, 
which are fixed by $\i_q$ for every $q\in U^{(1)}$.
Nevertheless we will now show that there is
a natural way to assign to them unique elements $q$.

\proclaim
{Lemma 5.15}
There exist unique elements $q_0\equiv 7(\mod 9)$ and
$q_1\equiv 4(\mod 9)$ in $U^{(1)}-U^{(2)}$ such that 
$\i_{q_0}$ and $\i_{q_1}$ have no nontrivial $3$-adic fixed points.
\endproclaim

\Pf
Fix $n\geq 2$.
By Theorem 5.8, the $q\in U^{(1)}-U^{(2)}$ for which $[\i_q]_{3^{2n-1}}$
has no rooted fixed points are exactly the $q$ for which $[\i_q]_{3^{2n-1}}$
fixes either $3^{n-1}$, if $q\equiv 7(\mod 9)$, or $1+3^{n-1}$,
if $q\equiv 4(\mod 9)$.  
By Theorem 5.8, together these are all of the $q$, with or without rooted fixed points,
for which $[\i_q]_{3^{2n-1}}$ has a fixed point $z$ with $v(z(z-1))=n-1$.
Since $n-1\leq (2n-1)-2=2n-3$, we may apply Proposition 5.2, which says these $q$
form the coset $q_n+3^n\Z_3$, where 
$q_n=1+a_1 3+\cdots+a_{n-1}3^{n-1}$,
$a_1=2$ if $q\equiv 7(\mod 9)$, $a_1=1$ if $q\equiv 4(\mod 9)$, 
and the $a_i\geq 2$ are uniquely determined elements of $\{0,1,2\}$.

It follows that
the elements $q\in U^{(1)}-U^{(2)}$ for which $[\i_q]_{3^{2n+1}}$ has
no rooted fixed points form the coset $q_{n+1}+3^{n+1}\Z_3$.
We claim that $q_{n+1}\equiv q_n(\mod 3^n)$.
If not then $[\i_{q_{n+1}}]_{3^{2n-1}}$ has a rooted fixed point, and
by Propagation Lemma 5.12, $[\i_{q_{n+1}}]_{3^{2n+1}}$ would have one too.
Thus $q_{n+1}=q_n+a_n 3^n$, for a unique $a_n\in\{0,1,2\}$.
In this way we define $q_n$ for all $n\geq 2$.
It is easily checked that
for all $r>n$ we have $q_r\equiv q_n(\mod 3^n)$,
so $v(q_r-q_n)\geq n$, and we have a well defined $3$-adic integer
$q=\lim_{n\to\infty}q_n$.
Since $[\i_{q_n}]_{3^{2n-1}}$ has no rooted fixed points,
neither does $[\i_{q_n}]_{3^n}$ by Propagation Lemma 5.12,
and since $q\equiv q_n(\mod 3^n)$, 
it follows that $[\i_q]_{3^n}$ has no fixed points, for all $n$.
Therefore $\i_q$ has no $3$-adic fixed points.
Since $a_1$ assumes the $2$ values $1$ and $2$,
and the $a_i$ are uniquely determined for $i\geq 2$, there are exactly 
$2$ of these $q$: $q_0\equiv 7(\mod 9)$ and $q_1\equiv 4(\mod 9)$.

\endpf

\flushpar
{\bf Remark 5.16.}
We can compute these numbers inductively.
By Theorem 5.8, the $[\i_q]_{3^{2n-1}}$ that have no rooted fixed points 
are exactly those that fix $z=3^{\lfloor\frac{2n-1}2\rfloor}=3^{n-1}$, 
if $q\equiv 7(\mod 9)$, or $z=1+3^{n-1}$, if $q\equiv 4(\mod 9)$.  Thus $v_0=v(z(z-1))={n-1}$,
and since $(2n-1)-v_0=n$, these $q$ 
form the coset $q_n+3^n\Z_3$, where $q_n=1+a_1 3+\cdots+a_{n-1}3^{n-1}$.
We then construct $q_n$ inductively for $n=2,3,\dots$.
For example, if $n=2$, it is easy to see that $q_2=1+2\cdot 3$,
i.e., $\i_7(3)\equiv 3(\mod 3^3)$.
To find $q_3$ we have to find the number $a_2\in\{0,1,2\}$ such that
for $q_3=q_2+a_2 3^2$, $[\i_{q_3}]_{3^5}$ fixes $3^2$,
and in general once $q_n$ has been found, there are only three candidates
to check to find $q_{n+1}$ such that $\i_{q_{n+1}}(3^{2n+1})$ fixes $3^{n-1}$.
Proceeding in this way, we compute
$$
\align
q_0&=1+ 2\cdot 3+ 3^2 + 3^6 + 2\cdot 3^7+\cdots
\\
q_1&=1 + 3 + 2\cdot 3^2 + 2\cdot 3^4 + 3^5 + 3^7 +\cdots
\endalign
$$
Note that the length of these computations grows exponentially.
For example, to find the $7$-th coefficient of $q_0$, we have to make $3^7=2187$
computations modulo $3^{15}=14,348,907$; the number $\i_q(7)$ already exceed $3^{15}$.

\

\flushpar
{\bf 6. Underlying Homeomorphisms and Isometries.}

\proclaim
{Theorem 6.1}
Suppose $q\in U^{(1)}-U^{(2)}$.
There exists a canonical homeomorphism
$$
\Psi:3\Z_3\cup(1+3\Z_3)\longrightarrow\;U^{(1)}-U^{(2)}
$$
defined
by $\Psi(0)=q_0$, $\Psi(1)=q_1$, where $q_0$ and $q_1$ are the elements of
Lemma 5.15, and by $\Psi(z)=q$ for all 
$z\in 3\Z_3\cup(1+3\Z_3)-\{0,1\}$, where $q\in U^{(1)}-U^{(2)}$
is the unique element such that $\i_q(z)=z$.
\endproclaim

\Pf
To show $\Psi$ is well defined it remains to show every element 
$z\in 3\Z_3\cup(1+3\Z_3)-\{0,1\}$ is a $3$-adic fixed point for 
some uniquely defined $q$.
By Theorem 5.13 it suffices to show $z$ is a rooted fixed point for
some $[\i_q]_{3^n}$, and since $1\leq v(z(z-1))<\infty$, this
is immediate by ``Existence of $q$'' Proposition 5.2, 
with $n>2\cdot v(z(z-1))+1$.
Since $q$ is uniquely defined, $\Psi$ is injective.

We show the image of $\Psi$ is dense in $U^{(1)}-U^{(2)}$.
To prove this it suffices to show that for every 
$q_r=1+a_1 3+\cdots+a_{r-1}3^{r-1}\in U^{(1)}-U^{(2)}$
with finite $3$-adic expansion, there exists a $q\in U^{(1)}-U^{(2)}$
in the image of $\Psi$, such that $q\equiv q_r(\mod 3^r)$.
For since any $q\in U^{(1)}-U^{(2)}$ may be approximated to arbitrary 
precision with a $q_r$ of finite $3$-adic expansion, this would show
any $q$ may be approximated to arbitrary precision with one that
has a $3$-adic fixed point.
Suppose by way of contradiction that $q_r$ is a counterexample,
i.e., $\i_q$ has no nontrivial $3$-adic fixed points
for any $q\in U^{(1)}-U^{(2)}$ such that $q\equiv q_r(\mod 3^r)$.
By Theorem 5.13, for these $q$, $[\i_q]_{3^n}$ has no rooted fixed points $z$
for any $n$.
But an easy count shows the number of $q$ of length $n\geq r$ extending $q_r$
is $3^{n-r}$, therefore we have at least $3^{n-r}$ elements $q\in U^{(1)}-U^{(2)}$
for which $[\i_q]_{3^n}$ has no rooted fixed points.
On the other hand, by Lemma 5.9, the number of generators $q\in U^{(1)}/U^{(n)}$ 
such that $q\equiv q_r(\mod 9)$
and $[\i_q]_{3^n}$ has no rooted fixed points is $3^{\lfloor\frac{n-1}2\rfloor}$.
For large enough $n$, this number is smaller than $3^{n-r}$; 
this occurs, for example, if $n>2r$.
Thus we have a contradiction, and we conclude the image of $\Psi$ is dense
in $U^{(1)}-U^{(2)}$.

Next we show that $\Psi$ is continuous.
First suppose $z\neq 0,1$ and $\Psi(z)=q$;
we'll show $\Psi$ is continuous at $z$.
Let $v_0=v(z(z-1))$, then $1\leq v_0<\infty$.
Fix any $N\geq v_0+2$.
Since the function $\i_q-\id$ is continuous, there exists a number
$n_0$ such that whenever $v(z-z')\geq n_0$ we have
$v(\i_q(z)-z-(\i_q(z')-z'))\geq N+v_0$.
Since $\i_q(z)=z$, this implies $\i_q(z')\equiv z'(\mod 3^{N+v_0})$,
hence $z$ and $z'$ are both modular fixed points for
$[\i_q]_{3^{N+v_0}}$.
If $q'=\Psi(z')$, then since $1\leq v_0<\infty$,
$q\equiv q'(\mod 3^N)$ by Proposition 5.2.
Therefore $v(z-z')\geq n_0$ implies $v(q-q')\geq N$,
and it follows that $\Psi$ is continuous at $z\neq 0,1$.

To show $\Psi$ is continuous at $0$ and $1$, suppose $z\in\{0,1\}$,
and let $q=\Psi(z)$.
Fix $N>>0$.
Suppose $z'\neq 0,1$, and $q'=\Psi(z')$.
If $v(z-z')\geq N-1$ then since $z$ equals $0$ or $1$ we have
$v(z'(z'-1))\geq N-1$,
and it follows by Theorem 5.13 that $[\i_{q'}]_{3^{2N-1}}$
has no rooted fixed points:
any rooted fixed points $z_0$ for $[\i_{q'}]_{3^{2N-1}}$ must satisfy
$v(z_0(z_0-1))<\frac{(2N-1)-1}2=N-1$ by definition, and by Theorem 5.13, $z'$
then satisfies $v(z'(z'-1))=v(z_0(z_0-1))$.
Since neither $[\i_q]_{3^{2N-1}}$ or $[\i_{q'}]_{3^{2N-1}}$
have rooted fixed points, they both fix $3^{\lfloor\frac{2N-1}2\rfloor}=3^{N-1}$,
if $q\equiv 7(\mod 9)$, or $1+3^{\lfloor\frac{2N-1}2\rfloor}=3^{N-1}$,
if $q\equiv 4(\mod 9)$.
By Proposition 5.2 we have $q\equiv q'(\mod 3^{(2N-1)-(N-1)})$,
i.e., $v(q-q')\geq N$.
Thus $v(z-z')\geq N-1$ implies $v(q-q')\geq N$.
Therefore $\Psi$ is continuous at $z\in\{0,1\}$.
We conclude that $\Psi$ is continuous.

Since $3\Z_3\cup(1+3\Z_3)$ is compact,
$\Psi$ is a closed map, in particular its image is closed.
Since its image is also dense, $\Psi$ is surjective. 
Since $\Psi$ is a continuous bijection that takes closed sets to closed sets, 
its inverse is continuous.
Therefore $\Psi$ is a homeomorphism.

\endpf

\proclaim
{Definition 6.2}
{\rm
Define 
$\Phi=\Psi^{-1}:U^{(1)}-U^{(2)}\;\longrightarrow\;3\Z_3\cup (1+3\Z_3)$,
the inverse of the map $\Psi$ of Theorem 6.1,
by assigning to each $q$ its $3$-adic fixed point.  
}
\endproclaim

\flushpar
{\bf Summary 6.3.}
The situation can be summarized as follows.  Any action of the additive group
$\Z_p$ on $\Z_p$ is defined by $1*1=q$ for some $q\in U^{(1)}$.
The canonical 1-cocycle $\i_q$ has a nontrivial $p$-adic fixed point if and only
if $p=3$, $q\in U^{(1)}-U^{(2)}$, and $q\not\in\{q_0,q_1\}$, the distinguished
$3$-adic integers of Lemma 5.15.
For these $q$ we have a unique $z_q\in\Z_3$, namely $\Phi(q)$.

We will explain how this fixed point governs the structure of the modular
fixed point sets for $[\i_q]_{3^n}$ for all $n$.
The fixed point $z_q$ satisfies $z_q\equiv 0$ or $1(\mod 3)$,
depending on whether $q\equiv 7$ or $4(\mod 9)$, respectively.
Thus we have $v_0:=v(z_q(z_q-1))\geq 1$.
The modular fixed point set for $[\i_q]_{3^n}$ is now given
explicitly by Theorem 5.8:  If $v_0\geq\frac{n-1}2$ and $q\equiv 7(\mod 9)$
then it is $3^{\lfloor\frac n2\rfloor}\Z_3\cup(1+3^{n-1}\Z_3)$.
In this case the $3$-adic fixed point $z_q$ does not manifest, in the sense that
there is no way to detect its $3$-value from the fixed point set, which has
the same structure as the modular fixed point set of any other 
$z_q$ with $v_0\geq\frac{n-1}2$.
If $v_0<\frac{n-1}2$ and $q\equiv 7(\mod 9)$ then $z_q$ appears as a ``rooted'' fixed point
$z_0$ that satisfies $v(z_0(z_0-1))=v_0$, and $z_q\equiv z_0(\mod 3^{n-v_0-1})$.
Moreover, $z_0$ appears with ``period'' $\tau=3^{n-v_0-1}$, meaning that $z_0+c\tau$
is a fixed point for $[\i_q]_{3^n}$ for all $c\in\Z_3$.
Additionally, the numbers $c\tau$ themselves are modular fixed points for 
$[\i_q]_{3^n}$, for all $c\in\Z_3$.
Aside from these fixed points, which are ``caused'' by the $3$-adic fixed point $z_q$,
there are the standard modular fixed points $1+3^{n-1}\Z_3$, which appear for
every prime $p$ and $q\in U^{(1)}$.
If $q\equiv 4(\mod 9)$, then the modular fixed point sets for various $3^n$ 
have the same broad characteristics as for the $q\equiv 7(\mod 9)$ case, except
that everything is ``shifted'' by $1$, as indicated in Theorem 5.8.

\proclaim
{Corollary 6.4}
Suppose $q\in U^{(1)}$, for general $p$.
Then the number of fixed points for $[\i_q]_{p^n}$ is asymptotically
stable for all $q$ and all $p$, except when $p=3$ and $q=q_0$ or $q=q_1$,
in which case the number grows without bound as $n$ goes to infinity.
\endproclaim

\Pf
We use ``asymptotically stable'' to mean the number is constant for large
enough $n$.
Corollary 5.10 shows that the number of fixed points is asymptotically 
stable when $p\neq 3$ or $q\in U^{(2)}$.
If $p=3$, $q\in U^{(1)}-U^{(2)}$, and $q\neq q_0,q_1$ then by
Theorem 6.1, $\i_q$ has a nontrivial $3$-adic fixed point $z_0$,
and for large enough $n$ this becomes a rooted fixed point for $[\i_q]_{3^n}$.
Explicitly, if $v_0=v(z_0(z_0-1))$, then as soon as $v_0<\frac{n-1}2$, 
$z_0$ is a rooted fixed point.
Thus the number of modular fixed points is asymptotically stable in this case.
When $q=q_0$ or $q=q_1$ there is no nontrivial $3$-adic fixed point,
hence no rooted fixed points for $[\i_q]_{3^n}$ for all $n$.
Case 4 of Corollary 5.10 then shows the number grows without bound
as $n$ goes to infinity.

\endpf

The most accessible values of $q$ are positive integers, and the most
accessible values of $\Z_3$ are rational numbers.  We have found
exactly one case where there is a rational $3$-adic fixed point for
$\i_q$ when $q$ is an integer.

\proclaim
{Proposition 6.5}
$\Phi(4)=-1/2$.
\endproclaim

\Pf
Since $4\in U^{(1)}$ we have $4^{-1/2}\in U^{(1)}$, and we see immediately that 
$4^{-1/2}=-1/2$.
Therefore $\i_4(-1/2)=(4^{-1/2}-1)/(4-1)=(-3/2)/3=-1/2$,
hence $\Phi(4)=-1/2$, as desired.

\endpf

We discover two isometries underlying the map 
$\Phi:U^{(1)}-U^{(2)}\to 3\Z_3\cup(1+3\Z_3)$.
The set $U^{(1)}-U^{(2)}$ is the union of the two cosets of the group $9\Z_3$,
$$
U^{(1)}-U^{(2)}=(4+9\Z_3)\cup(7+9\Z_3),
$$
and by Corollary 5.4,
$\Phi$ takes $4+9\Z_3$ onto $1+3\Z_3$, and $7+9\Z_3$ onto $3\Z_3$.
Let $f_1:\Z_3\isim 4+9\Z_3$,
$g_1:\Z_3\isim 7+9\Z_3$, $f_2:1+3\Z_3\isim\Z_3$, and
$g_2:3\Z_3\isim\Z_3$ be the obvious
topological isomorphisms. 
Define 
$$
F:\Z_3\longrightarrow\Z_3,\qquad G:\Z_3\longrightarrow\Z_3
$$ 
to be the compositions $F=f_2\cdot\Phi\cdot f_1$ and $G=g_2\cdot\Phi\cdot g_1$.

\proclaim
{Theorem 6.6}
The functions $F$ and $G$ are isometries.
\endproclaim

\Pf
We'll show $v(F(x)-F(x'))=v(x-x')$ for all $x,x'\in\Z_3$; 
the proof for $G$ is similar.
Choose $x,x'\in\Z_3$ such that,
in the above notation, $q=f_1(x)$ and $q'=f_1(x')$ avoid the distinguished
element $q_0$.
We have seen that $z_q=\Phi(q)$ and $z_{q'}=\Phi(q')$ are elements of $1+3\Z_3$.
Set $v_0=v(z_q-1)$ and $v_0'=v(z_{q'}-1)$, then $1\leq v_0,v_0'<\infty$.
Let $n>2v_0+1$.
Then $1\leq v_0\leq n-2$, and by Proposition 5.2,
if $q\equiv q'(\mod 3^{n-v_0})$ then $\i_{q'}(z_q)\equiv z_q(\mod 3^n)$.
Since $n>2v_0+1$, by Definition 5.6, $z_q$ is a rooted fixed point of 
$[\i_{q'}]_{3^n}$.  By Theorem 5.13, $v_0=v_0'$, and
$z_{q'}$ is also a rooted fixed point of $[\i_{q'}]_{3^n}$.
By Theorem 5.8, $z_q\equiv z_{q'}(\mod 3^{n-v_0-1})$.
Conversely, if $z_q\equiv z_{q'}(\mod 3^{n-v_0-1})$, then by Theorem 5.8, 
$z_{q'}$ is a rooted fixed point of $\i_q$,
and by Theorem 5.13, $v_0=v_0'$.
Since $z_{q'}$ is a fixed point of $\i_{q'}$, we have
$q\equiv q'(\mod 3^{n-v_0'})$ by Proposition 5.2, hence $q\equiv q'(\mod 3^{n-v_0})$.

Thus $q\equiv q'(\mod 3^{n-v_0})$ if and only if $z_q\equiv z_{q'}(\mod 3^{n-v_0-1})$.
Equivalently,
$\frac 19(q-4)\equiv \frac 19(q'-4)(\mod 3^{n-v_0-2})$ if and only if $\frac 13(z_q-1)
\equiv \frac 13(z_{q'}-1)(\mod 3^{n-v_0-2})$.
That is, by definition, for all $m\geq v_0$, 
$x\equiv x'(\mod 3^m)$ if and only if $F(x)\equiv F(x')(\mod 3^m)$.
Thus $v(x-x')=v(F(x)-F(x'))$.

We have shown that $F$ is an isometry on the punctured disk $\Z_3-\{x_0\}$,
where $x_0=f_1^{-1}(q_0)$.
By definition $F$ is a homeomorphism on the whole disk, as the composition of
the homeomorphisms $f_1,\Phi$, and $f_2$.
It follows by the continuity of the $3$-adic metric that $F$ is an isometry on 
all of $\Z_3$.
For if $x_0=\lim_{n\to\infty}s_n$, where $s_n\in\Z_3-\{x_0\}$,
then $v(x-x_0)=\lim_{n\to\infty}v(x-s_n)$,
hence $v(F(x)-F(x_0))=\lim_{n\to\infty}v(F(x)-F(s_n))$ since $F$ is an isometry
on the punctured disk, hence $v(F(x)-F(x_0))=v(F(x)-F(x_0))$, as desired.
This completes the proof.

\endpf

\

\flushpar
{\bf 7. Examples.}

We work some examples for $p=3$ and $q=4$.
We have noted already that the $3$-adic fixed point for $\i_4$ is 
$\Phi(4)=-1/2=1+3+3^2+\cdots$.
By Theorem 4.2, we always have the modular fixed point pairs $c 3^{n-1}$ and $1+c 3^{n-1}$.
The $3$-adic fixed point $-1/2$
determines all remaining modular fixed points as follows.
The associated value is $v_0=v(-1/2(-3/2))=1$, so we expect to see rooted fixed
points for $\i_4(\mod 3^n)$ for $n>2v_0+1=3$, i.e., starting with $\i_4(\mod 3^4)$.
The period for $3^n$ is $3^{n-v_0-1}=3^{n-2}$,
so since $4\equiv 4(\mod 9)$, by Theorem 5.8, $1+c 3^{n-2}$ is fixed for all $c\in\N$, 
and then there will be rooted fixed points determined
by the residue of $-1/2(\mod 3^{n-2})$, added to the multiples of $3^{n-2}$.

Here is the sequence of values $\i_4(z)(\mod 3^4)$ 
from $z=0$ to $z=3^4+1=82$.
$$
\matrix
\boxed 0 & \boxed 1 & 5 & 21 & \boxed 4 & 17 & 69 & 34 & 
56 & 63 & \boxed{10} & 41 & 3 & \boxed{13} & 53 & 51 & 43 & 11\\
45 & \boxed{19} & 77 & 66 & \boxed{22} & 8 & 33 & 52 & 47 & \boxed{27} &
\boxed{28} & 32 & 48 & \boxed{31} & 44 & 15 & 61 & 2\\
9 & \boxed{37} & 68 & 30 & \boxed{40} & 80 & 78 & 70 & 38 & 72 & \boxed{46} &
23 & 12 & \boxed{49} & 35 & 60 & 79 & 74\\
\boxed{54} & \boxed{55} & 59 & 75 & \boxed{58} & 71 & 42 & 7 & 29 & 36 &
\boxed{64} & 14 & 57 & \boxed{67} & 26 & 24 & 16 & 65\\
18 & \boxed{73} & 50 & 39 & \boxed{76} & 62 & 6 & 25 & 20 & \boxed{0} & \boxed{1}
 & ...\\
\endmatrix
$$
The fixed point pairs are $c 3^3$ and $1+c3^3$, for $c=0,1,2$.
The rooted fixed points are $\{4,13,22,31,\dots\}$,
and the period is $3^{4-v_0-1}=3^2=9$. 
Note there are exactly $2\cdot 3^{v_0+1}+3=21$ fixed points modulo $3^4$, 
as required by Corollary 5.10.
To contrast, we list the values of $\i_4(\mod 3^5)$. 
To detect the pattern we only need list $\i_4(z)(\mod 3^5)$ from $z=0$ to
$z$ equals $1$ plus the period, i.e.,
$z=1+3^{5-v_0-1}=28$. 
$$
\matrix
\boxed{\ssize 0} & \boxed{\ssize 1} &\ssize 5 &\ssize 21 &\ssize 85 &\ssize 98 &\ssize 150 
&\ssize 115 &\ssize 
218 &\ssize 144 &\ssize 91 &\ssize 122 &\ssize 3 & \boxed{\ssize 13} &\ssize 53 &\ssize 213 &\ssize 124 &\ssize 11\\
\ssize 45 &\ssize 181 &\ssize 239 &\ssize 228 &\ssize 184 &\ssize 8 &\ssize 33 &\ssize 133 &\ssize 47 &\ssize 189 &\ssize
\boxed{\ssize 28} &\ssize 113 &\ssize 210 &\ssize 112 &\ssize 206 &\ssize 96 &\ssize 142 &\ssize 83 \\ 
\ssize 90 &\ssize 118 &\ssize 230 &\ssize 192 &\ssize \boxed{\ssize 40} &\ssize 161 &\ssize 159 &\ssize 151 &\ssize 119 &\ssize 234 &\ssize 
208 &\ssize 104 &\ssize 174 &\ssize 211 &\ssize 116 &\ssize 222 &\ssize 160 &\ssize 155 \\
\ssize 135 & \boxed{\ssize 55} &\ssize 221 &\ssize 156 &\ssize 189 &\ssize 71 &\ssize 42 &\ssize 169 &\ssize 191 &\ssize 36 &\ssize
145 &\ssize 95 &\ssize 138 & \boxed{\ssize 67} &\ssize 26 &\ssize 105 &\ssize 178 &\ssize 227 \\
\ssize 180 &\ssize 235 &\ssize
212 &\ssize 120 &\ssize 238 &\ssize 224 &\ssize 168 &\ssize 187 &\ssize 20 &\ssize \boxed{\ssize 81} &\ssize \boxed{\ssize 82} &\ssize
86 &\ssize 102 &\ssize 166 &\ssize 179 &\ssize 231 &\ssize 196 &\ssize 56 \\ \ssize 225 
&\ssize
172 &\ssize 203 &\ssize 84 &\ssize \boxed{\ssize 94} &\ssize 134 &\ssize 51 &\ssize 205 &\ssize 92 &\ssize 126 &\ssize 19 &\ssize 77 &\ssize 66 &\ssize 
22 &\ssize 89 &\ssize 114 &\ssize 214 &\ssize 128 \\ \ssize 27 &\ssize
\boxed{\ssize 109} &\ssize 194 &\ssize 43 &\ssize 193 &\ssize 44 &\ssize 177 &\ssize 223 &\ssize 164 &\ssize 171 &\ssize 199 &\ssize 68 &\ssize 30 &
\boxed{\ssize 121} &\ssize 242 &\ssize 240 &\ssize 232 &\ssize 200 \\
\ssize 72 &\ssize 46 &\ssize 185 &\ssize 12 &\ssize 49 &\ssize 197 &\ssize 60 &\ssize
241 &\ssize 236 &\ssize 216 & \boxed{\ssize 136} &\ssize 59 &\ssize 237 &
\dots \\
\endmatrix
$$
In addition to the fixed point pairs,
the other fixed points will be $27+13=40$, $1+2\cdot 27=55$, $2\cdot 27+13=67$, etc.
Here are the values of $\i_4(\mod 3^6)$, up to $1$ plus the period $3^4$.
$$
\matrix
\boxed{\ssize 0} & \ssize \boxed{\ssize 1}& \ssize 5 & \ssize 21 & \ssize 85 & 
\ssize 341 & \ssize 636 & \ssize 358 & \ssize 
704 & \ssize 630 & \ssize 334 & \ssize 608 & \ssize 246 & \ssize 256 & \ssize 296 & 
\ssize 456 & \ssize 367 & \ssize 11\\
\ssize 45 & \ssize 181 & \ssize 725 & \ssize 714 & \ssize 670 & \ssize 494 & \ssize 519 & 
\ssize 619 & \ssize 290 & \ssize 432 & \ssize
271 & \ssize 356 & \ssize 696 & \ssize 598 & \ssize 206 & \ssize 96 & \ssize 385 & 
\ssize 83\\
\ssize 333 & \ssize 604 & \ssize 230 & \ssize 192 & \ssize \boxed{\ssize 40} & \ssize 161 & 
\ssize 645 & \ssize 394 & \ssize 119 & \ssize 477 & \ssize 451 & \ssize
347 & \ssize 660 & \ssize 454 & \ssize 359 & \ssize 708 & \ssize 646 & \ssize 398\\
\ssize 135 & \ssize 541 & \ssize 707 & \ssize 642 & \ssize 382 & \ssize 71 & \ssize 285 & 
\ssize 412 & \ssize 191 & \ssize 36 & \ssize
145 & \ssize 581 & \ssize 138 & \ssize 553 & \ssize 26 & \ssize 105 & \ssize 421 & 
\ssize 227\\
\ssize 180 & \ssize 721 & \ssize 698 & \ssize 606 & \ssize 238 & \ssize 224 & \ssize 168 & 
\ssize 673 & \ssize 506 & \ssize 567 & \ssize \boxed{\ssize 82}
 & \dots\\
\endmatrix
$$

\newpage

\flushpar
{\bf References.}

\roster
\item"{\cite{A}}"
Arens, R.:
{\it Homeomorphisms preserving measure in a group},
Ann. of Math., {\bf 60}, no. 3, (1954), pp. 454--457.
\item"{\cite{B}}"
Bishop, E.:
{\it Isometries of the $p$-adic numbers},
J. Ramanujan Math. Soc. {\bf 8} (1993), no. 1-2, 1--5.
\item"{\cite{C}}"
Conrad, K.:
{\it A $q$-analogue of Mahler expansions. I,}
Adv. Math.  {\bf 153}  (2000),  no. 2, 185--230.
\item"{\cite{D}}"
Dieudonn\'e, J.: 
{\it  Sur les fonctions continues $p$-adiques},
Bull. Sci. Math. (2) {\bf 68} (1944), 79--95
\item"{\cite{F}}"
Fray, R. D.:
{\it Congruence properties of ordinary and $q$-binomial coefficients,}
Duke Math. J. {\bf 34} (1967), 467--480.
\item"{\cite{J}}"
Jackson, F. H.:
{\it $q$-difference equations},
Amer. J. Math. {\bf 32} (1910), 305--314.
\item"{\cite{G}}"
Granville, A.:  {\it http://www.dms.umontreal.ca/$\sim$andrew/Binomial/elementary.html}
\item"{\cite{M}}"
Mahler, K.:
{\it An interpolation series for continuous functions of a $p$-adic variable},
J. Reine Angew. Math. {\bf 199} (1958) 23--34.
\item"{\cite{N}}"
Neukirch, J.:
{\it Algebraic Number Theory,} 
Grundlehren der Mathematischen Wissenschaften, 322, 
Springer-Verlag, Berlin, 1999.
\item"{\cite{S}}"
Serre, J.-P.:
{\it A Course in Arithmetic}, Springer-Verlag, New York, 1973.
\item"{\cite{Su}}"
 Sushchanski\u i, V. I.:
{\it Standard subgroups of the isometry group of the metric space of $p$-adic integers},
V\B isnik Ki\"iv. Un\B iiv. Ser. Mat. Mekh. 
{\bf 117} no. 30 (1988), 100--107.
\endroster

\bye